\documentclass[reqno]{amsart}
\usepackage{color}
\usepackage{comment}
\usepackage{pdfpages}
\usepackage{graphicx}
\usepackage{dcolumn}
\usepackage{dsfont}

\RequirePackage[numbers]{natbib}
\RequirePackage[colorlinks=true, pdfstartview=FitV, linkcolor=blue,
  citecolor=blue, urlcolor=blue]{hyperref}
\usepackage{paralist}

\usepackage{tikz} 
\usepackage{dcolumn}
\usetikzlibrary{arrows,positioning}
\tikzset{
    >=stealth',
    pil/.style={
           ->,
           thick,
           shorten <=2pt,
           shorten >=2pt,}
}

\usepackage{graphicx}
\graphicspath{{./figures/}}
\usepackage{amsfonts}
\usepackage{amsmath}
\usepackage{amsthm}
\usepackage{amssymb}
\usepackage{amsbsy}
\usepackage{epsfig}
\usepackage{natbib, mathrsfs}
\usepackage{verbatim}
\usepackage[latin1]{inputenc}
\usepackage{mhequ}
\usepackage{algorithm}
\usepackage{algorithmic}

\usepackage{import}
\usepackage{enumitem} 
\usepackage{caption,subcaption}

      \let\e=\varepsilon

  \let\s=\sigma

\renewcommand{\leq}{\;\leqslant\;}                   
\renewcommand{\geq}{\;\geqslant\;}                   
\newcommand{\sumtwo}[2]{\sum_{\substack{#1 \\ #2}}} 

\newcommand{\var}{\operatorname{Var}}

\newcommand\sind[1]{^{(#1)}}
\newcommand\weakconv{\stackrel{d}{\rightarrow}}
\newcommand\sumw{\mathcal{Z}}

\newcommand\tdlim{N/L \rightarrow \rho}

\newcommand\sparseapmts{\ref{aspt_A1_pos_limit_exists} - \ref{aspt_A3_theta} }

\newcommand{\N}{\mathbb N}

\newcommand{\bbN}{{\ensuremath{\mathbb N}} }

\newcommand{\Z}{{\ensuremath{\mathbb Z}} }

\definecolor{WowColor}{rgb}{.75,0,.75}
\definecolor{SubtleColor}{rgb}{0.9,0,0}

\newcounter{margincounter}

\newcounter{latercounter}

\newtheorem{theorem}{Theorem}[section]
\newtheorem{lemma}[theorem]{Lemma}
\newtheorem{proposition}[theorem]{Proposition}

\newtheorem{remark}[theorem]{Remark}

\newtheorem{definition}[theorem]{Definition}

\newtheorem*{question*}{Question}

\newtheorem*{remark*}{Remark}
\newtheorem*{idefinition*}{Definition}
\newtheorem*{example*}{Example}

\begin{document}
\title[Poisson-Dirichlet asymptotics in condensing particle systems]{Poisson-Dirichlet asymptotics in\\ condensing particle systems} 
\author[P. Chleboun]{Paul Chleboun}
\author[S. Gabriel]{Simon Gabriel}
\author[S. Grosskinsky]{Stefan Grosskinsky}
\address{P. Chleboun, Department of Statistics, University of Warwick,  Coventry,  CV4 7AL,  United Kingdom } 
\email{paul.i.chleboun@warwick.ac.uk}
\address{S. Gabriel, Mathematics Institute, University of Warwick, Coventry CV4 7AL, United Kingdom } 
\email{simon.gabriel@warwick.ac.uk}
\address{S. Grosskinsky, Department of Mathematics, University of Augsburg, 86159 Augsburg, Germany} 
\email{stefan.grosskinsky@math.uni-augsburg.de}
\thanks{S. Gabriel acknowledges financial support from EPSRC through grant EP/R513374/1. S. Grosskinsky is grateful to Technical University of Delft, where part of this research was carried out.}
\keywords{Poisson-Dirichlet distribution, split-merge dynamics, random partitions, equivalence of ensembles, interacting particle systems, condensation}

\begin{abstract}
We study measures on random partitions, arising from condensing stochastic particle systems with stationary product distributions. We provide fairly general conditions on the stationary weights, which lead to
Poisson-Dirichlet statistics of the condensed phase in the thermodynamic limit. The Poisson-Dirichlet distribution is known to be the unique reversible measure of split-merge dynamics for random partitions, which we use to characterize the limit law. 
We also establish concentration results for the macroscopic phase, using size-biased sampling techniques and the equivalence of ensembles to characterize the bulk distribution of the system. 


\end{abstract}

\maketitle


\setcounter{tocdepth}{1}
\tableofcontents

\section{Introduction and results}
\subsection{Mathematical setting and motivation}

The results presented in this paper are motivated by the study of interacting particle systems.  We consider finite systems consisting of $N$ particles on $L$ sites indexed by the set $\Lambda$, where $|\Lambda| = L$. For simplicity we take $\Lambda = \{1,\ldots,L\}$ in the remainder of this paper. 
The space of such particle configurations $\eta =(\eta_1 ,\ldots ,\eta_L )$ with $\eta_x \in \N_0$ is given by 
\begin{equation*}
\Omega_{L,N}:= \left\{ 
\eta \in \mathbb{N}_0^L: \sum_{x=1}^L \eta_x = N
\right\}\, ,
\end{equation*}
which we equip with the discrete topology.
The dynamics of systems we consider are assumed to be irreducible Markov processes on $\Omega_{L,N}$, conserving only the quantities $N$ and $L$. Thus, there exists a unique invariant distribution $\pi_{L,N}$ of the system on $\Omega_{L,N}$, which is called the \emph{canonical distribution}.

We will focus on models where the {canonical distributions} are of product form 
\begin{equation}\label{wpi}
\pi_{L,N}[d\eta] = 
\begin{cases}
\displaystyle\frac{1}{Z_{L,N}} \prod_{x=1}^L w_L(\eta_x)\, d\eta\  & \text{if } \sum_{x=1}^L \eta_x = N\, ,\\
0 & \text{otherwise}.
\end{cases} 
\end{equation}
Here $d\eta$ denotes the counting measure on $\Omega_{L,N}$ and $(w_L(n))_{n\in\mathbb{N}_0}$ is a sequence of positive weights, possibly depending on the system size $L$. The normalising constant, called \emph{canonical partition function}, is given as
\[
Z_{L,N} =\sum_{\eta\in\Omega_{L,N}} \prod_{x=1}^L w_L(\eta_x)\, .
\]
Note that the weights $(w_L(n))_{n\in\mathbb{N}_0}$ are independent of the site $x$, thus the $\pi_{L,N}$ are permutation invariant and in particular spatially homogeneous, so that single-site marginals $\pi_{L,N}[\eta_x\in \cdot]$ do not depend on $x$.\\

We are primarily interested in the limiting behaviour of $\pi_{L,N}$ in the thermodynamic limit $N,L\to \infty $ such that $\tfrac{N}{L}$ converges to $\rho\geq 0$, which we will subsequently abbreviate by $\tdlim$. 
Assume for now that the weak limit of the single-site marginals exists for all $\rho\geq 0$, 
\begin{align}\label{eq_Piweakconv}
    \pi_{L,N} [\eta_x \in \cdot \,]\to \nu_{\rho}[\,\cdot\,] \quad \text{as }\tdlim\, ,
\end{align}
and the limit is a probability measure on $\N_0$. 
This implies in particular convergence of the expectations for bounded functions $f: \mathbb{N}_0\mapsto \mathbb{R}$, i.e.
\begin{align*}
\sum_{n=0}^{\infty} f(n) \pi_{L,N} [\eta_x=n ]
\to
\nu_{\rho}(f)\,  ,   
\end{align*}
where we write $\nu_{\rho}(f)$ for the expectation of $f$ under $\nu_{\rho}$. Sometimes we will also write $\nu_{\rho}(f(\eta_x))$ for the corresponding expectation if needed for clarity. 
Looking at a single site's expected occupation number under $\pi_{L,N}$, we see that due to spatial homogeneity we have 
$$
\pi_{L,N}(\eta_x)
=
\frac{1}{L} \sum_{y=1}^L \pi_{L,N}(\eta_y)
=
\frac{1}{L}  \pi_{L,N}\Big(\sum_{y=1}^L\eta_y\Big)
=\frac{N}{L}\to\rho\, ,
$$ 
when taking the thermodynamic limit. 
However, because the identity $f(n)=n$ is an unbounded function,
we cannot guarantee that the particle density of the system is conserved in the limit $\nu_{\rho}$ and $\nu_\rho (\eta_x )$ may be strictly smaller than $\rho$.
This phenomenon is known as condensation. 
\begin{definition}[Condensation]\label{def_cond}
A system characterised by spatially homogeneous canonical distributions $(\pi_{L,N})_{L,N}$  exhibits \textbf{condensation} in the thermodynamic limit $\tdlim$ if $\nu_\rho$ in \eqref{eq_Piweakconv} exists and
\begin{equation*}
\nu_{\rho}(\eta_x ) 
< \rho = \lim_{\tdlim} \pi_{L,N} (\eta_x)\, .
\end{equation*}
Furthermore, we say that the system has a \textbf{condensation transition} with \textbf{critical density} $\rho_c\geq 0$ if 
\begin{equation*}
\nu_{\rho}(\eta_x ) 
\begin{cases}
= \rho & \quad \text{ if } \rho < \rho_c\, , \\
< \rho & \quad \text{ if } \rho>\rho_c\, .
\end{cases}
\end{equation*}
\end{definition}
In the context of stochastic particle systems, condensation means that a positive fraction of the total density $\rho$ is not observed in the thermodynamic limit, since it concentrates on sites with diverging occupation numbers called the \emph{condensed phase}. Clearly, the number of such sites has a vanishing volume fraction and does not contribute to the weak limit $\nu_{\rho}$, which describes the distribution of the \emph{background} or \emph{bulk phase}. 

Condensation in homogeneous stochastic particle systems has been studied previously in great generality, partially reviewed e.g.\ in \cite{ChGr13,EW14,G19}. Early results are formulated in the context of zero-range processes in \cite{DrGoCa98,Ev00,Go03} and in \cite{JeMaPi00,GrScSp03,ArLo08,ArGrLo13} on a rigorous level, where the condensed phase concentrates on a single lattice site. Our goal here is to understand details of the condensed phase when it extends over more than one site and exhibits a non-trivial structure. Such structures have previously been observed as a result of spatial correlations \cite{WaEv12,WaSoJa09,ThTaCaBl10} and as a result of $L$-dependent stationary weights, with a soft cut-off for site occupation numbers under zero-range dynamics \cite{ScEvMu08} or in the inclusion process \cite{JCG19}.\\

On the level of particle configurations the condensed phase disappears in the thermodynamic limit due to its vanishing volume fraction. To study its structure, it is more useful to interpret a configuration as an ordered partition of the total mass. 
For models of type \eqref{wpi} partitions and particle configurations are equivalent, since $\pi_{L,N}$ is permutation invariant and the underlying lattice structure is irrelevant.


We will represent particle configurations rescaled by the total mass $N$ as ordered partitions of the unit interval $[0,1]$
on the set
\begin{align}\label{eq_def_nabla_bar}
\overline{\nabla}
:= 
\left\{p=
(p_i )_i \in [0,1]^{\mathbb{N}}
:
\sum_{i=1}^{\infty}p_i \leq 1 
\quad \text{and}\quad
p_1\geq p_2 \geq \cdots 
\right\}\, .
\end{align}
We use the map $T=T(L,N):\Omega_{L,N} \to \overline{\nabla}$ with
\begin{equation}\label{eq_def_T}
T(\eta ):=\frac{1}{N}\left( \widehat{\eta}_1,\ldots ,\widehat{\eta}_L ,0,\ldots\right)\, ,
\end{equation}
where $\widehat{\eta} =(\widehat{\eta}_1,\ldots, \widehat{\eta}_L)$ denotes the entries in $\eta$ in decreasing order with $\widehat{\eta}_1\geq \widehat{\eta}_2\geq \ldots \geq \widehat{\eta}_L$. Since any permutation of entries in $\eta$ yields the same partition in $\overline{\nabla}$, the map $T$ is not injective. Thus, the push-forward measure of $\pi_{L,N}$ under $T$ on $\overline{\nabla}$ is given by 
\begin{equation}\label{eq_def_mu}
\mu_{L,N}[dp] := \pi_{L,N} \circ T^{-1} [dp]
= \pi_{L,N} [d(Np)] \,|T^{-1}(\{ p \})|\, ,
\end{equation}
with $Np$ denoting any configuration in $\Omega_{L,N}$ inducing the finite ordered partition $p\in\overline{\nabla}$. 
Note that $\mu_{L,N}$ concentrates on finite partitions with at most $L$ non-zero entries and $T^{-1} (\{ p\} )=\emptyset$ otherwise. 
In fact, the $\mu_{L,N}$ further concentrate on the subset where $\sum_{i=1}^\infty p_i=1$. However, this space is not compact, unlike $\overline{\nabla}$ which is compact w.r.t.\ the product topology by Tychonoff's theorem, ensuring existence of subsequential weak limits of $\mu_{L,N}$ in the thermodynamic limit $N/L\to\rho$. 
The objective of this article is to identify general assumptions on the weights $(w_L)_L$, such that 
for $\rho>0$ large enough $\mu_{L,N}$ converges weakly to a Poisson-Dirichlet distribution as $\tdlim$. Details on this distribution are introduced in Section \ref{sec_math_setting}. \\

The starting point of our analysis is the recent paper 
\cite{JCG19} in which weights of the form 
\begin{equation}\label{eq_ip_weights}
w_L(n)=\frac{\Gamma(n+d)}{n! \Gamma(d)}\, ,
\end{equation}
with $d=d(L)\in\mathbb{R}$  such that $\lim_{L\to \infty}dL=\theta \in (0,\infty)$,
are considered. Such weights emerge for example from the dynamics of the inclusion process introduced in \cite{GKR07,CGGR13}, which can also be applied in population genetics as a multi-species Moran model \cite{Mo58}, with the above scaling of the parameter $d$ corresponding to a small mutation rate. 
With weights \eqref{eq_ip_weights} the system \eqref{wpi} exhibits a condensation transition with critical density $\rho_c =0$, leaving an empty bulk behind, i.e. $\nu_{\rho}(\eta_x) = 0$. For the condensed phase in this model we have
\begin{align}\label{eq_thm1_JCG}
    \mu_{L,N}\weakconv \text{PD}(\theta)\, , \quad \text{ as }\,\tdlim\,\text{ for all }\rho >0\, ,
\end{align}
where PD$(\theta)$ denotes the Poisson-Dirichlet distribution with parameter $\theta$ \cite[Theorem 1]{JCG19}. 
The proof uses the fact that $\pi_{L,N}$ with weights \eqref{eq_ip_weights} is a Dirichlet multinomial distribution which permits an exact, simple expression for the corresponding partition function $Z_{L,N}$. This leads to exact expressions for the distribution of size-biased marginals which characterize the Poisson-Dirichlet limit (cf.\  Section \ref{sec_PD}). 

Our main result provides a generalization to models with more general weights that do not lead to exact expressions for $Z_{L,N}$, and with non-trivial bulk distribution where $0<\nu_{\rho}(\eta_x)<\rho$.
%
%
In our proof, we not only make use of the Poisson-Dirichlet distribution's characterisation via size-biased sampling, which was essential for the arguments in \cite{JCG19}, but also use the characterisation as the unique reversible distribution under split-merge dynamics as explained in Section \ref{sec_PD}. 
This allows us to avoid explicit expressions or approximations of the partition function $Z_{L,N}$ which are not always at hand. 
Our approach is motivated by a recent paper by Ioffe and T\'{o}th \cite{IoTo20}, where the embedding of integer configurations into partitions of $[0,1]$ was used to show convergence of cycle-length processes of stationary random stirring to the split-merge dynamics.

\subsection{Main results}
We recall from \eqref{wpi} that the $\pi_{L,N}$ are probability measures on $\Omega_{L,N}$ given by 
\begin{align*}
    \pi_{L,N}[d\eta] = 
\frac{1}{Z_{L,N}} \prod_{x=1}^L w_L(\eta_x)\, d\eta\, .
\end{align*}
For our first result we fix a density $\rho >0$ and choose $N,L\to\infty$ such that $\tdlim$.
We assume that
\medskip
\begin{enumerate}[label=(A\arabic*)]
 \item\label{aspt_A1_pos_limit_exists} $w_L (n)\geq 0$ for all $n\in \mathbb{N}$
and the limit 
\begin{equation*}
\lim_{L\to\infty} w_L(n)=:w(n)\geq 0\quad\mbox{exists for all fixed  }n\in\mathbb{N}_0\, ,
\end{equation*}
such that $w$ is summable and non-trivial,\\
\end{enumerate}
and a weak form of the equivalence of ensembles:
\medskip
\begin{enumerate}[label=(A\arabic*),resume]
\item\label{aspt_A2_equiv_ensemble} 
The limiting probability distribution \eqref{eq_Piweakconv} exists and is of the form

\begin{equation*}
 \pi_{L,N}[\eta_1=n]\to \nu_\rho [n]=
 \frac{w(n)\phi^n}{\sumw}\qquad\forall\,n\geq 0\, ,
\end{equation*}
for some $\phi> 0$, and $\sumw =\sum_n w(n)\phi^n \in (0,\infty )$ is the corresponding normalising constant.\\
\end{enumerate}

\begin{remark}
The equivalence of ensembles is the main mathematical framework to understand the large scale behaviour of statistical mechanics models, and in particular to show condensation as in Definition \ref{def_cond} (see Section \ref{sec_proof_specific_result} for details). Assumption \ref{aspt_A2_equiv_ensemble} has therefore been established for all homogeneous particle systems that are known to exhibit condensation (see citations above). 
We will see in our second result Theorem \ref{THEO_SPECIALISED_MAIN}, that equivalence of ensembles and Assumption \ref{aspt_A2_equiv_ensemble} can be shown for a large class of models under slightly stronger assumptions on convergence of the stationary weights.
\end{remark}

\begin{remark}
It is natural to ask if $\nu_\rho$, cf. \eqref{eq_Piweakconv},  can always be described by an exponential change of measure w.r.t. $w$, as is assumed in \ref{aspt_A2_equiv_ensemble}. 
%
%
One possible sufficient condition to recover the postulated form in \ref{aspt_A2_equiv_ensemble} is to require that
$\nu_\rho$ and the measure $w$ defined by the limiting weights in \ref{aspt_A1_pos_limit_exists} are equivalent in the sense of measures, i.e. $\nu_\rho$ and $ w$ have the same support,
 and $w(0),\, w(1)>0$. Then we can define
$$
\phi := \frac{\nu_{\rho}[1]}{w(1)}\frac{w(0)}{\nu_{\rho}[0]}>0
\quad \text{and}\quad
\sumw = \frac{w(0)}{\nu_{\rho}[0]}\in (0,\infty )
$$
and use a telescopic product argument for the ratio of partition functions to see that
\begin{align*}
    \frac{Z_{L-1,N-n}}{Z_{L,N}}=
    \frac{Z_{L-1,N}}{Z_{L,N}}\prod_{k=1}^n \frac{Z_{L-1,N-k}}{Z_{L-1,N-k+1}}\to \frac{1}{\sumw}\phi^n\, .
\end{align*}
Thus, we recover Assumption \ref{aspt_A2_equiv_ensemble} in this case, and a similar argument works in the degenerate case $w(0)>0$ and $w(k)=0$ for all $k\geq 1$, which applies for the inclusion process with weights \eqref{eq_ip_weights}. 
We believe that this connection between $\nu_\rho$ and $w$ holds under more general conditions, but a more detailed discussion is out of the scope of this paper.
\end{remark}

In order to formulate our main result with simple notation, we assume without loss of generality that $\phi= \sumw=1$, since we can absorb $\phi$ and $\sumw$ into the weights by 
\begin{equation*}
\widetilde{w}_L(n) := \frac{w_L(n) \phi^n}{\sumw} \quad \text{and }
\quad 
\widetilde{w}(n) := \frac{w(n) \phi^n}{\sumw}\,.
\end{equation*}
So $w$ can be assumed to be the probability mass function of $\nu_\rho$, i.e.
\begin{equation}\label{eq_nu_rho_equal_w}
    \nu_\rho [\eta_x =n]=w(n)\quad\mbox{for all }n\geq 0\,  .
\end{equation}

We are interested in the macroscopic part of the condensed phase, i.e. the distribution of occupation numbers $\eta_x$ that scale linearly with the total mass $N$ in the system when taking the thermodynamic limit. 
The structure of this macroscopic phase will depend on the asymptotic behaviour of the stationary weights. In order to see Poisson-Dirichlet statistics in this phase, the weights $w_L(n)$ must scale like $(nL)^{-1}$, for at least all $n$ which are visible under macroscopic rescaling:
\medskip
\begin{enumerate}[label=(A\arabic*),resume]
\item \label{aspt_A3_theta}
 there exists $\theta> 0$ such that for all $\varepsilon\in (0,1)$ we have 
\begin{equation*}
\sup_{\varepsilon N \leq n \leq N}
\left|
n w_L(n)L-\theta
\right|\to 0 \qquad \text{as }\tdlim\,  .
\end{equation*}
\end{enumerate}
\medskip
For the second part of Theorem \ref{THEO_MAIN} we also impose a second moment condition:
\medskip
\begin{enumerate}[label=(A\arabic*),resume]
\item \label{aspt_A4_concentration}
The limit 
\begin{equation}\label{eq_def_alpha}
\alpha^2 := \frac{1+\theta}{\rho} \lim_{\tdlim} \frac{\pi_{L,N}(\eta_x^2)}{N}\in [0,1]\quad\text{exists}\, .
\end{equation}
\end{enumerate}
\medskip
In Section \ref{sec_split_merge} we will see that
\begin{equation*}
    \alpha = \mu (\| p\|_1 )\, ,\quad\mbox{with}\quad \| p\|_1 =\sum_{j=1}^\infty p_j\, ,
\end{equation*}
i.e. $\alpha$ coincides with the expected total mass fraction for each accumulation point $\mu$ of the measures $(\mu_{L,N})_{L,N}$, and that the variance of $\| p\|_1$ vanishes. Therefore, Assumption \ref{aspt_A4_concentration} guarantees that the macroscopic phase is well defined in the thermodynamic limit, excluding fluctuations of mass towards other scales, and plays an important role when identifying the accumulation points $\mu$. This leads to our first main result.



\begin{theorem}\label{THEO_MAIN}
Let $\rho> 0$ and $(w_L)_{L}$ be a sequence of weights satisfying
\sparseapmts for some $\theta\in (0,1]$ and let $\mu$ be an accumulation point of 
the laws of mass partitions $(\mu_{L,N})_{L,N}$ defined in \eqref{eq_def_mu}. Then
\[
\mu =\text{PD}_{[0,\alpha_\mu ]}(\theta) \quad\mbox{with}\quad \alpha_\mu =\mu\big( \| p\|_1 \big) \in [0,1]\ ,
\]
where PD$_{[0,\alpha' ]}(\theta)$ denotes the Poisson-Dirichlet distribution with parameter $\theta$, concentrating on partitions of the interval $[0,\alpha ']$, which depends on the accumulation point. 
If in addition \ref{aspt_A4_concentration} holds then
\begin{equation*}
(\mu_{L,N})_{L,N}
\weakconv \text{PD}_{[0,\alpha]}(\theta)\, , \quad \text{ as }\quad\tdlim \, .
\end{equation*}

\end{theorem}

\begin{remark}
(a) If $\theta =0$, assumption \ref{aspt_A3_theta} does not specify the leading order limiting behaviour of the weights. Our proof can cover this case, if we assume in addition that $n\mapsto nw_L (n)$ is a regularly varying function\footnote{i.e.  $nw_L (n) /\big(\lambda nw_L ([\lambda n])\big)\to C\in (0,\infty )$ for all $\lambda >0$ as $n\to\infty$} for large enough $L$ (see \eqref{eq_punchline} in the proof). This is consistent with choosing fixed weights of the form $w_L(n)=w(n)=n^{-b} /\sumw$ for $b>2$ whenever $n>0$ and $w_L(0)=1$. 
Indeed, for this choice, as part of a larger class of sub-exponential weights, it is a well known result that the condensed phase consists of a single cluster \cite{Ev00,GrScSp03,ArLo08}, which can be interpreted as a degenerate Poisson-Dirichlet distribution with $\theta=0$.

(b) The result also covers the case $\alpha=0$ where the macroscopic phase is empty and its limiting distribution is trivial. This includes models that do not condense at all, or where the condensed phase concentrates on sub-macroscopic scales such as for certain models with spatial correlations \cite{WaEv12,WaSoJa09,ThTaCaBl10}. Of course our result does not say anything interesting in this case, since there is no mass on the macroscopic scale.

(c) All results in this paper extend to arbitrary $\theta> 0$ under the assumption that PD$(\theta)$ is the unique reversible distribution for the split-merge dynamics introduced in Section \ref{sec_math_setting}. It seems widely accepted that this is indeed the case, though to the authors' best knowledge no proof exists for $\theta>1$.
\end{remark}

Assumption \ref{aspt_A3_theta} is the core premise which guarantees the Poisson-Dirichlet limit of the macroscopic phase, and is consistent with the scaling of weights \eqref{eq_ip_weights} for the inclusion process
$$
w_L(n) =\frac{\Gamma(n+d)}{n!\Gamma(d)}\simeq \frac{\theta}{n\, L}\,,\quad\text{if we set}\quad d=\frac{\theta}{L}\,  .
$$
This scaling implies that a size-biased sample of a macroscopic cluster has the stationary weights $nw_L (n)\simeq \theta /L$ which are independent of $n$. That means, picking a particle uniformly at random, the size of its cluster is uniformly distributed. This is the trademark of Poisson-Dirichlet statistics, and the distribution can only be normalized in the scaling limit due to the factor $1/L$. 
Thus, to get a non-trivial macroscopic phase with $\theta >0$, it is necessary that the weights $w_L$ depend on the system size $L$. 
The role of $\theta$ and more details on size-biased sampling will be given in Section \ref{sec_math_setting}.

\begin{remark}
A system where \sparseapmts are satisfied but \ref{aspt_A4_concentration} is difficult to verify, is for example given by weights of the form 
\begin{align}\label{eq_intermediate_weights}
    w_L(n) = 
    w(n)
    \mathds{1}\{n\leq  N^{1/2}\}
    +\frac{\theta}{n\, L}\mathds{1}\{ n> N^{1/2}\}
    \,,
\end{align}
with $w$ being an arbitrary probability mass function on $\N$ (not necessarily of finite support).
In this case, to determine if \ref{aspt_A4_concentration} holds,
careful evaluation of the  second moment condition  would be required that would depend on the choice of $w(n)$. 
Theorem~\ref{THEO_MAIN} still allows to characterise the macroscopic part of the condensate as Poisson-Dirichlet, with the catch that it is possibly trivial with $\alpha =0$.
\end{remark}

Recall that in Theorem \ref{THEO_MAIN} we have fixed the density $\rho >0$ and it provides a very general result that also includes trivial cases without condensation. But we had to assume the equivalence of ensembles in \ref{aspt_A2_equiv_ensemble} and regularity of the macroscopic phase in \ref{aspt_A4_concentration}, which are not easy to check in general (if not established already for particular models). Strengthening the requirements \ref{aspt_A1_pos_limit_exists} and \ref{aspt_A3_theta} on the stationary weights $(w_L)_L$, we can use Theorem \ref{THEO_MAIN} to show a stronger but more specialized result, including the equivalence of ensembles and regularity of the macroscopic phase in the conclusion.

\begin{theorem}\label{THEO_SPECIALISED_MAIN}
Assume $(w_L)_L$ is a sequence of non-negative weights satisfying the following two conditions:
\begin{enumerate}[label=(B\arabic*)]
 \item \label{aspt_B1_uniform_conv}
 $(w_L)_L$ converges in the sup-norm, $\|\cdot\|_\infty$, to a sequence $w$ such that $$\sum_{n=0}^{\infty}w(n)=1$$
 and either $w(0)=1$ or 
 \begin{align}\label{eq:bern_ass}
 w(0)>0\quad\mbox{and}\quad\sup_{n}[ w(n-1)\wedge w(n)]>0\,.
 \end{align}

\item \label{aspt_B3_theta} There exists some $\theta \in (0,1]$ such that
\begin{equation*}
\lim_{J\to \infty }\lim_{L\to\infty}
\sup_{n > J}\left|n w_L(n) L -\theta \right| =  0\,  .
\end{equation*}
\end{enumerate}
Then the system exhibits a condensation transition according to Definition \ref{def_cond} with critical density
\[\displaystyle\rho_c := \sum_{n=0}^{\infty}n w(n)\in [0,\infty )\,.\]
Furthermore,  we have bulk density $\nu_\rho (\eta_x )=\rho_c$ for all $\rho\geq\rho_c$ and
\begin{align*}
(\mu_{L,N})_{L,N}
\weakconv \text{PD}_{[0,\alpha]}(\theta)\, , \quad \text{ as }\quad\tdlim\geq\rho_c \, ,
\end{align*}
with $\alpha=\alpha(\rho) = 1-\frac{\rho_c}{\rho}$.
\end{theorem}

\begin{remark}\label{def_A}
By assumption \ref{aspt_B3_theta}, there exists an $A\in \mathbb{N}_0$  such that
$$
\lim_{L\to\infty}
\sup_{ n > A}\left|n w_L(n) L -\theta \right| \leq  C\quad\text{with an arbitrary constant }C>0\, .
$$
So for each $n \geq A$ we have $w_L(n) \to 0$ as $L\to \infty$, and hence by assumption \ref{aspt_B1_uniform_conv}
 \begin{equation}\label{wlim}
     w(n)=0 \quad \text{for all } n >A\,,
 \end{equation}
and the limiting distribution $w$ can only have finite support.
In this sense, Theorem~\ref{THEO_MAIN} is more general because it allows for arbitrary limiting distributions of possibly infinite support, at the cost of loosing control over intermediate scales. Recall \eqref{eq_intermediate_weights} for an example.
\end{remark}

Clearly, the restriction in \ref{aspt_B1_uniform_conv} that $w$ can be interpreted as a probability mass function is for notational convenience, we could just assume summability. 
Furthermore, with $\theta >0$ all models covered by this result have a macroscopic phase with non-trivial structure, excluding systems where the latter concentrates on a single site which have been studied previously (see citations above).
Condition \eqref{eq:bern_ass} is necessary to avoid lattice effects and establish the equivalence of ensembles for the bulk part of the distribution (see Proposition \ref{prop_single_site_marg}).
In the special case that $w(0) = 1$ there is a simpler proof of the equivalence of ensembles result.


We want to stress that single-site-condensation in models is typically due to a strong enough attraction between particles, whereas for systems covered in Theorem \ref{THEO_SPECIALISED_MAIN} particle attraction alone is too weak, and condensation only occurs in combination with particle expulsion from the bulk as represented by condition \eqref{wlim}  on the limiting weights. 
This strict exclusion condition for occupation numbers larger than $A$ in the limiting weights prevents clustering of particles on sub-macroscopic scales. 
It should be possible to weaken this, but some form of bulk exclusion is essential for condensation with non-trivial macroscopic phase in models with stationary product measures. In Section \ref{sec_application} we provide an intuitive explanation of this in terms of dynamics of generic particle systems covered by our result.


\subsection{Key steps of the proofs}

The essential steps of the proof of Theorem \ref{THEO_MAIN} may be summarised as follows.
By compactness we know that $(\mu_{L,N})_{L,N}$ has weak accumulation points. In order to determine the limit points' distributions, we prove that $(\mu_{L,N})_{L,N}$ is approximately reversible w.r.t. a discrete split-merge dynamics. These discrete dynamics converge to the generator of the coagulation-fragmentation process with split-merge dynamics, which we will introduce in Section \ref{sec_math_setting}, see \eqref{eq_discrete_gen_split_merge}.
Lastly, we use size-biased sampling together with a disintegration argument to prove that the corresponding limit points concentrate and therefore have a Poisson-Dirichlet law. This is due to the fact that the Poisson-Dirichlet distribution is the unique distribution which concentrates and is reversible w.r.t. the limiting split-merge dynamics mentioned above \cite{DMWZZ04,Sc05}. Here, we say that $\mu$ concentrates if $\|p\|_1=\alpha$ $\mu$-a.s. for some $\alpha\in [0,1]$.\\

Essentially, Theorem \ref{THEO_SPECIALISED_MAIN} is a direct application of Theorem \ref{THEO_MAIN}. Additionally, the stronger assumptions allow us to establish the equivalence of ensembles in Appendix \ref{sec:equivensembles}, which in our case implies the condensation transition.
The proof is based on the application of a local central limit theorem (LCLT) which, together with a relative entropy bound, shows convergence of single-site marginals to a distribution independent of the particle density $\rho\geq \rho_c$. \\

The remainder of the paper will be structured as follows: in Section~\ref{sec_math_setting} we will give a short review on Poisson-Dirichlet distributions and size-biased sampling. Section~\ref{sec_split_merge} will focus on the proof of Theorem~\ref{THEO_MAIN} whereas in Section \ref{sec_proof_specific_result} 
we state the proof of Theorem~\ref{THEO_SPECIALISED_MAIN}. 
Lastly, we discuss possible applications of Theorem \ref{THEO_SPECIALISED_MAIN} to a family of zero-range and generalized inclusion processes in Section \ref{sec_application}. 
In Appendix \ref{sec:equivensembles} we give a brief introduction to grand-canonical ensembles before proving equivalence of ensembles for size-dependent weights under a sub-exponential growth condition.
\\
\medskip

\section{Background on partitions}\label{sec_math_setting}

\subsection{The Poisson-Dirichlet distribution}\label{sec_PD}

The Poisson-Dirichlet (PD) distribution is a one-parameter family of probability measures on the space of ordered partitions $\nabla =\nabla_{[0,1]}$ of the unit interval, where for any $\alpha >0$ we denote
\begin{equation*}
{\nabla}_{[0,\alpha]}
:= 
\left\{p=
(p_i)_i \in [0,1]^{\mathbb{N}}
:
\| p\|_1  = \alpha
\quad \text{and}\quad
p_1\geq p_2 \geq \cdots 
\right\}\,  .
\end{equation*}
Note that elements $p$ in $\nabla$ are not partitions themselves but induce partitions of the form $\{ [0,p_1),[p_1,p_1+p_2),\ldots \}$.
The family of measures was first introduced by Kingman \cite{Ki75} in the study of random distributions on countably infinite sets, motivated by Bayesian inference and decision theory. 
Apart from the original construction as a limit of Dirichlet distributions,
the PD distribution can be more intuitively constructed via a stick-breaking procedure. Let $U_1,U_2,\ldots$ be independent Beta$(1,\theta)$-distributed random variables and define
\begin{equation*}
V_1 := U_1\,,\; \;
V_2 := (1-U_1)U_2\,,\;\;
V_3 := (1-U_1)(1-U_2)U_3\,,
\ldots,
\end{equation*}
i.e. we start with a stick of unit length and continue by breaking a random fraction of $U_1$ apart. Then we do the same with the remaining part of the stick and iterate.
The resulting random vector $V=(V_i)_{i\geq 1}$ is said to be GEM($\theta$)-distributed, named after Griffiths \cite{Gr80,Gr88} Engen \cite{En78} and McCloskey \cite{McC65}. Reordering the entries of $V$ in decreasing order yields $\widehat{V}$ which is known to be PD($\theta$)-distributed (see e.g.\ \cite{Fe10}).

Note the two special cases, $\theta =1$ where the $U_i$'s are uniformly distributed on the interval $[0,1]$, and $\theta = 0$ where the $U_i$'s are degenerated point-measures on one and hence $\widehat{V}=V=(1,0,0,\ldots)$. 
Clearly, the choice of the interval $[0,1]$ is arbitrary and one can construct PD and GEM distributions on intervals $[0,\alpha]$ for arbitrary $\alpha>0$ just by rescaling
\begin{equation*}
    p\sim PD_{[0,\alpha]}(\theta)\quad\Leftrightarrow\quad p/\alpha :=(p_1 /\alpha ,p_2 /\alpha ,\ldots )\sim PD_{[0,1]}(\theta)\,  ,
\end{equation*}
and analogously for  $GEM_{[0,\alpha]}(\theta)$. 
Since its introduction in \cite{Ki75} the PD distribution emerged first in population biology \cite{Ki75,EtKu81}, before appearing in statistical mechanics  \cite{KMRTSZ07,GLU12,IoTo20} and interacting particle systems \cite{JCG19}. In particular, the statistics of cycle length distributions in spatial permutations can be linked to condensation phenomena in quantum-mechanical models, see e.g. \cite{GUW11,BeUe11,zeindler} and references therein. While those models often involve weights that satisfy a decay condition similar to \ref{aspt_A3_theta}, they consider a different scaling limit with only one diverging parameter, and our results are not directly applicable.\\

Besides its characterisation via the GEM-construction, the PD distribution was furthermore found to be the unique invariant (in fact reversible) measure on $\nabla$ of the coagulation-fragmentation process with split-merge dynamics for $\theta \in (0,1]$. 
This is a Markov process on the state space $\overline{\nabla}$
with infinitesimal generator given by 
\begin{equation}\label{eq_discrete_gen_split_merge}
\mathcal{G}_{\theta} f(p) = 
\sum_{i\neq j}p_i p_j \left[ f(\widehat{M}_{ij}p) -f(p) \right]
+
\theta
\sum_i p_i^2 \left[  \int_0^1 f(\widehat{S}_{i}^{u}p) du   -f(p) \right]\, .
\end{equation}
Here $\widehat{M}_{ij}p$ denotes the operator that merges the parts $p_i$ and $p_j$ to a single block of size $p_i+p_j$ and then reorders the partition to maintain the decreasing order. On the other hand, $\widehat{S}_{i}^{u}p$ defines the operation of splitting $p_i$ into two blocks of size $up_i$ and $(1-u)p_i$ before reordering the resulting partition.

\begin{proposition}[\cite{DMWZZ04,Sc05}]\label{prop_PD_unique_inv}
For $\theta \in [0,1]$, the Poisson-Dirichlet distribution PD($\theta$) is the unique invariant measure on $\nabla$ with respect to 
split-merge dynamics defined by $\mathcal{G}_{\theta}$, and it is also reversible.
\end{proposition}

Since the generator in \eqref{eq_discrete_gen_split_merge} conserves the total mass of partitions,
it is clear that 
there exist stationary distributions for split-merge dynamics on $\nabla_{[0,\alpha ]}$ for all $\alpha >0$, which are unique and equal to $PD_{[0,\alpha ]} (\theta )$ with the above result. 
In general, the set $C_b (\overline{\nabla})$ of bounded continuous functions is the natural domain for the (Feller) Markov semigroup associated to split-merge dynamics. Under the product topology on $\overline{\nabla}$, cf. \eqref{eq_def_nabla_bar}, these include in particular bounded cylinder functions, which depend only on finitely many entries of a partition, and for all such functions the generator \eqref{eq_discrete_gen_split_merge} is well defined (see also \cite[Lemma 4]{MWZZ02}). 

Originally, the split-merge process was constructed in discrete time, see \cite{MWZZ02}, the extension to continuous time can be found in \cite[Section 7.4]{GUW11}.
The uniqueness of the invariant measure when $\theta=1$ was proven in \cite{DMWZZ04} by Zerner, Zeitouni, Mayer-Wolf and Diaconis. An alternative technique allowed Schramm to extend this result to $\theta\in (0,1]$, see \cite{Sc05} and Theorem 7.1 in \cite{GUW11}. Lastly, consider the case $\theta=0$, then clearly $\delta_{(1,0,\ldots)}$ is invariant because $\mathcal{G}_{0}$ only consists of the merge term and there cannot exist another invariant measure on $\nabla$.

\subsection{Size-biased sampling}\label{sec_size_biased}


Partitions in $\nabla$ can be interpreted as probability mass functions themselves, which allows for a natural size-biased resampling of its elements.
Given $p= (p_i)_{i\in\mathbb{N}}\in \nabla$ we sample an index $i\in \mathbb{N}$ at random according to $(p_i)_{i\in\mathbb{N}}$. Continuing this procedure, while renormalising the remaining partition to a total mass of one in each round, we construct a so-called \emph{size-biased sample} $\widetilde{p}$ of $p$. For given $p\in\nabla$, $\widetilde{p}$ is a random element of the unordered set
\begin{equation*}
{\Delta}
:= 
\left\{p=
(p_i)_i \in [0,1]^{\mathbb{N}}
:
\| p\|_1  = 1
\right\} \, ,
\end{equation*}
and we denote its distribution by $\sigma_p$. Of course $\widetilde p$ can also be defined in the same way for $p\in\Delta$.

The above procedure can be generalised to partitions in the compact space
\begin{equation*}
\overline{\Delta}:= 
\left\{p=
(p_i)_i \in [0,1]^{\mathbb{N}}
:
\| p\|_1  \leq 1
\right\}\,  ,
\end{equation*}
which includes in particular $\overline{\nabla} = \bigcup_{\alpha\in [0,1]}\nabla_{[0,\alpha]}$ 
that we already introduced in \eqref{eq_def_nabla_bar}.
More precisely, fix an element $p\in\overline{\Delta}$. 
In the following, $q=\widetilde p$ will denote the size-biased sample with distribution $\sigma_p$ on $\overline{\Delta}$ which is defined recursively:
\begin{itemize}
\item the first entry of $q$ is assigned the value\footnote{
We refrain from using an equal sign instead of '$\leftarrow$', since this could lead to mathematically wrong statements. For example, if $p=(\tfrac{1}{2},\tfrac{1}{2},0,\ldots)$ then $q_1=\tfrac{1}{2}$ with probability $\frac12 +\frac12 =1$ which is reflected by our notation in \eqref{eq_q_first_sample}, but would read $q_1=\tfrac{1}{2}$ with probability $\tfrac{1}{2}$ when replacing '$\leftarrow$' with '$=$'.
In a fully rigorous construction of size-biased samples we actually sample the index $j$ at random and not the value $p_j$, see \cite{Gn98} for the full construction and more details.
Because our analysis does not differentiate between entries of the same size, we omit this step to significantly simplify notation and assign the value $p_j$ directly. 
} 
\begin{equation}\label{eq_q_first_sample}
q_1 \leftarrow
\begin{cases}
p_j &\quad \text{ w.p. } p_j \text{ for all }j\in\mathbb{N}\, , \\
0  &\quad \text{ w.p. } 1-\|p\|_1\, ,
\end{cases}
\end{equation}
\item for $i>1$, let $V$ be the set of indices of $p$ assigned to $q_k$ for $k<i$, then 
\begin{equation}\label{eq_q_later_samples}
q_i
 \leftarrow
\begin{cases}
p_j &\quad \text{ w.p. } \displaystyle\frac{p_j}{1-\sum_{k=1}^{i-1}q_k} \text{ for all }j\in \mathbb{N}\setminus V\, , \\
&\\
0  &\quad \text{ w.p. } \displaystyle\frac{1-\|p\|_1}{1-\sum_{k=1}^{i-1}q_k}\, .
\end{cases}
\end{equation}
\end{itemize}
In the case where $p$ consists of finitely many non-zero components only, we sample zeros in each iteration after exhausting all non-zero components. In contrast to size-biased sampling on $\Delta$ with $\| p\|_1 =1$, 
which is usually defined in terms of shuffling indices of the original sequence, 
size-biased sampling of  $p\in\overline{\Delta}$ assumes a non-exhaustive reservoir of zeros from which we pick with probability proportional to $1-\|p\|_1$ in each round.
For given $p\in\overline{\Delta}$, the distribution $\sigma_p [dq]$ then denotes the law on $\overline{\Delta}$ of $q$ defined above, and concentrates on partitions with $\| q\|_1 = \| p\|_1$.



For an arbitrary probability measure $\mu \in \mathcal{M}_1(\overline{\nabla})$ we then define its size-biased distribution $\sigma(\mu)$ as the law 
\begin{equation}\label{eq_def_sigma_mu}
\sigma(\mu)[dq] := \int_{\overline{\nabla}} \sigma_p[dq] \mu [dp]\, .
\end{equation}
One interesting result regarding size-biased distributions with reservoirs of zeros, is that weak convergence of measures on $\overline{\nabla}$ implies weak convergence of the corresponding size-biased distributions.

\begin{lemma}\label{lem_size_biased_conv}
If a sequence $(\mu_n)_{n\in\mathbb{N}}$ of probability measures on $\overline{\nabla}$ converges weakly to a measure $\mu$, then also
$\sigma(\mu_n)\weakconv \sigma(\mu)$ on $\overline{\Delta}$.
\end{lemma}

Originally this result was stated in \cite[Theorem 1]{DoJo89} with a flawed construction and proof. A correct
proof of the lemma can be found in \cite[Theorem 1]{Gn98}  along with a nice exposition on size-biased sampling.\\



From the stick-breaking construction of the PD distribution it is easy to see (e.g.\ in \cite{Fe10}) that the size-biased distribution of PD$(\theta)$ on $\nabla_{[0,1]}$ is precisely the GEM($\theta$) distribution. Considering on the other hand PD$_{[0,\alpha]}(\theta)$ on $\nabla_{[0,\alpha]}$, its size-biased distribution, as defined in \eqref{eq_q_first_sample} and \eqref{eq_q_later_samples}, contains $0$-elements and does not coincide with GEM$_{[0,\alpha ]} (\theta )$ whenever $\alpha<1$. 
%
However, this connection still holds for a 
modified (positive) size-biasing $\tilde p'$ without $0$-elements, defined again via scaling.  
For $p\in\overline{\Delta}\setminus\{ 0\}$ (i.e.\ we have $\| p\|_1 >0$), the positive size-biased sample is defined as
\begin{equation}\label{eq_def_pos_sb}
    \tilde p' :=\| p\|_1 q \quad\mbox{where}\quad q\sim \sigma_{p/\| p\|_1} \,  .
\end{equation}
Therefore the law of $\tilde p'$ on $\overline{\Delta}\setminus\{ 0\}$ is $\sigma'_p [dq] = \sigma_{p/\| p\|_1} \circ\frac{1}{\| p\|_1} [dq]$, the pushforward measure under rescaling of the size-biasing $\sigma_{p/\| p\|_1}$ on $\Delta$, which does not contain any $0$-elements. As for $\sigma_p$, $\sigma'_p [dq]$ concentrates on partitions with $\| q\|_1 =\| p\|_1$, and for a distribution $\mu$ on $\overline{\Delta}$ we define $\sigma' (\mu)$ analogously to \eqref{eq_def_sigma_mu}. 


Note that $\sigma(\mu) = \sigma'(\mu)$ on $\mathcal{M}_1(\nabla)$ since $\|p\|_1=1$ $\mu$-a.s., and in this case we will write ${\mu}[\widetilde{p}\in\cdot]$ instead to ease notation.
Furthermore, we can recover the finite dimensional marginals of $\sigma'(\mu)$ from $\sigma(\mu)$ by ignoring zeros sampled w.r.t. the latter measure \cite{Gn98}. We will only use this fact for the first marginal and state its derivation for completeness: for fixed $p\in \overline{\nabla}\setminus \{0\}$
\begin{align}\label{eq_sb_fm}
    \sigma_p[q_1 \in \cdot \;| \; q_1>0 ]
    &=
    \frac{\sigma_p[q_1 \in \cdot\,,\  q_1>0 ]}{\sigma_p[ q_1>0 ]}
    =
    \sum_{j=1}^\infty \frac{p_j \mathds{1}\{p_j>0,\, p_j\in \cdot\, \} }{\|p\|_1}\\
    &=
    \sum_{j=1}^\infty \frac{p_j }{\|p\|_1} \mathds{1}\{ p_j\in \cdot\,\}
    =
    \sigma_{p/\|p\|_1}[q_1\in \tfrac{\cdot}{\|p\|_1} ]
    =
    \sigma_p'[q_1\in \cdot\,]
\nonumber\,,
\end{align}
where we started with the definition of $\sigma_p$, cf. \eqref{eq_q_first_sample}, and rewrote the conditional probability in terms of $\sigma_p'$.
The statement for general $\mu$ follows by averaging. 

Returning to the case of the PD distribution $\mu \sim PD_{[0,\alpha]}(\theta)$, we can retrieve the first marginal of the corresponding GEM distribution on $[0,\alpha]$:
\begin{equation*}
\sigma'(\mu)[{q}_1\in\cdot]=
\sigma(\mu)\left[q_1 \in\cdot\;\middle|\; q_1>0 \right] = \mathrm{Beta}_{[0,\alpha]}(1,\theta)\,  .
\end{equation*}
Note also that for positive size-biased distributions the equivalent statement of Lemma \ref{lem_size_biased_conv} does not hold, because loss of mass of the corresponding sequences $p\sim \mu_n$ may occur. This would correspond to a positive probability of sampling a block size of zero in the limit, which is precisely the case for distributions $\mu_{L,N}$ we study in this paper.
\\


We will work with size-biased sampling not only on $\overline{\Delta}$ but also on $\Omega_{L,N}$, which corresponds to uniformly picking a particle and sampling the occupation number on its site. 
Using our definition from above (cf.\ \cite[Section 2.3]{JCG19} for a more detailed construction), the size-biased distribution of $\pi_{L,N}$ is given by 
\begin{equation*}
\sigma(\pi_{L,N})[d\eta]
:= 
\sigma(\mu_{L,N})\left[d \left(\frac{\eta}{N} \right) \right]\,  ,
\end{equation*}
where $\mu_{L,N}$ \eqref{eq_def_mu} is the distribution of the ordered partition corresponding to a particle configuration $\eta$. 
Because $\mu_{L,N}\in\mathcal{M}_1(\nabla)$, we will write $\pi_{L,N}[\widetilde{\eta}\in \cdot]$ instead of $\sigma(\pi_{L,N})$, and note that due to the product structure of \eqref{wpi} we have e.g.\ for the first marginal
\begin{equation}\label{sbform}
    \pi_{L,N} [\widetilde{\eta}_1 =n] =\frac{L}{N} n\,w_L (n) \frac{Z_{L-1,N-n}}{Z_{L,N}}\quad\mbox{for all }n=0,\ldots ,N\,  .
\end{equation}
Here, we used the notation $\widetilde{\eta}$ for the size-biased configuration, in the same manner as for partitions on $\overline{\Delta}$.

\medskip

\section{Proof of Theorem \ref{THEO_MAIN}}\label{sec_split_merge}
Recall the canonical distributions $\pi_{L,N}$ on $\Omega_{L,N}$ given by
\begin{equation*}
\pi_{L,N}[d\eta] 
=
\frac{1}{Z_{L,N}}\prod_{z=1}^L w_L(\eta_z) d\eta
\end{equation*}
and their macroscopic counterparts $\mu_{L,N}$ on $\nabla\subset \overline{\nabla}$, see \eqref{eq_def_mu}.
The proof of Theorem \ref{THEO_MAIN} can be broken down into two main steps. We begin by showing that every accumulation point of $(\mu_{L,N})_{L,N}$ is reversible with respect to the infinitesimal generator $\mathcal{G}_{\theta}$ of the 
split-merge process \eqref{eq_discrete_gen_split_merge}. The second step consists of 
proving that each limiting measure on $\overline{\nabla}$ in fact concentrates on $\nabla_{[0,\alpha]}$ for some $\alpha$. Together with Proposition \ref{prop_PD_unique_inv} this implies that the limiting measure must be a Poisson-Dirichlet distribution on the interval $[0,\alpha]$.\\

\subsection{Reversibility of weak accumulation points}


The space $\overline{\nabla}$ of ordered \mbox{(sub-)}partitions is compact w.r.t. the product topology on $[0,1]^{\mathbb{N}}$, which implies compactness of the space  $\mathcal{M}_1(\overline{\nabla})$ w.r.t. the topology induced by weak convergence. Therefore, every subsequence of $(\mu_{L,N})_{L,N}$ has a weakly convergent subsequence and $(\mu_{L,N})_{L,N}$ has at least one accumulation point.\\

We prepare the proof of Theorem \ref{THEO_MAIN} which is given at the end of this Section by the following intermediate result.
\begin{proposition}\label{prop_pd_limit}
Consider weights $(w_L)_L$ satisfying assumptions \sparseapmts with $\theta> 0$. 
All weak accumulation points of $(\mu_{L,N})_{L,N}$ are reversible w.r.t. $\mathcal{G}_{\theta}$.
\end{proposition}

In the remainder of this subsection, we set the stage for the proof of the proposition, which can be found after the statement of Lemma~\ref{lem_split_merge_transf_A}.\\

Because $\mu_{L,N}$ concentrates on discrete partitions 
and by Assumption \ref{aspt_A3_theta} we only control the weights on the macroscopic scale, 
it will be convenient to consider a discrete version $\mathcal{G}\sind{N,\varepsilon}_{\theta} :C_b (\overline{\nabla})\to C_b (\overline{\nabla})$ of the split-merge process corresponding to $\mathcal{G}_{\theta}$ \eqref{eq_discrete_gen_split_merge}, which only acts on the parts of the partition exceeding a fixed size of $\varepsilon\in(0,1)$,
\begin{equation*}
\begin{split}
\mathcal{G}\sind{N,\varepsilon}_{\theta} f(p)
&:= 
\frac{N}{N-1}\sum_{i\neq j}p_i p_j
\mathds{1}\{ p_i,p_j \geq \varepsilon \}
 \left[ f(\widehat{M}_{ij} p) -f(p) \right]\\
&\qquad +
\frac{\theta}{N-1}\sum_i p_i  
\mathds{1}\{ p_i \geq 2\varepsilon \}
  \sum_{k=\varepsilon N }^{N(p_i -\varepsilon)}\left[ f(\widehat{S}_{i}^{\frac{k}{Np_i}}p)   -f(p) \right]\, .
  \end{split}
\end{equation*}
As we will see Lemma \ref{lem_conv_generators}, it suffices to control the weights on macroscopic scales larger than $\varepsilon$ in order to observe a PD limit.

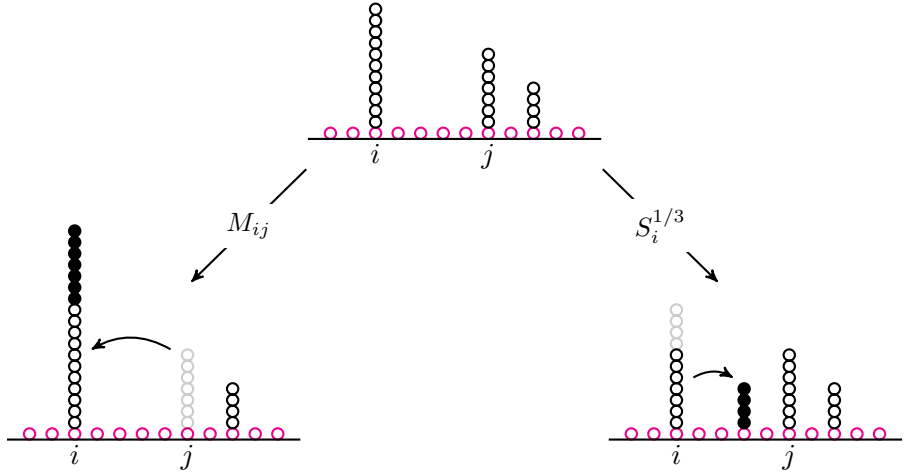
\begin{figure}[H]
    \centering 
\begin{subfigure}{0.95\textwidth}
  \begin{tikzpicture}
\node[inner sep=0pt, anchor= south] (config) at (4,4)
    {
	\begin{tikzpicture}[scale = 0.3]
\draw (3,-.20) node[anchor=north]  {\large{$i$}}
(8,-.20) node[anchor=north]  {\large{$j$}}
;

\draw[thick, magenta] (1,0.25) circle (0.25) {} 
(2,0.25) circle (0.25) {}
(3,0.25) circle (0.25) {}
(4,0.25) circle (0.25) {}
(5,0.25) circle (0.25) {}
(6,0.25) circle (0.25) {}
(7,0.25) circle (0.25) {}
(8,0.25) circle (0.25) {}
(9,0.25) circle (0.25) {}
(10,0.25) circle (0.25) {}
(11,0.25) circle (0.25) {}
(12,0.25) circle (0.25) {};

\draw[thick] (3,0.75) circle (0.25) {} 
(3,1.25) circle (0.25) {} 
(3,1.75) circle (0.25) {} 
(3,2.25) circle (0.25) {} 
(3,2.75) circle (0.25) {} 
(3,3.25) circle (0.25) {} 
(3,3.75) circle (0.25) {} 
(3,4.25) circle (0.25) {} 
(3,4.75) circle (0.25) {} 
(3,5.25) circle (0.25) {} 
(3,5.75) circle (0.25) {} ;

\draw[thick] (8,0.75) circle (0.25) {} 
(8,1.25) circle (0.25) {} 
(8,1.75) circle (0.25) {} 
(8,2.25) circle (0.25) {} 
(8,2.75) circle (0.25) {} 
(8,3.25) circle (0.25) {} 
(8,3.75) circle (0.25) {} ;

\draw[thick] (10,0.75) circle (0.25) {} 
(10,1.25) circle (0.25) {} 
(10,1.75) circle (0.25) {} 
(10,2.25) circle (0.25) {} ;

\draw[thick] (0,0) -- (13,0);

\end{tikzpicture}    
    };
\node[inner sep=0pt, anchor= south] (configmerge) at (0,0)
    {
    \begin{tikzpicture}[scale = 0.3]
\draw (3,-.20) node[anchor=north]  {\large{$i$}}
(8,-.20) node[anchor=north]  {\large{$j$}}
;
\draw[thick, magenta] (1,0.25) circle (0.25) {} 
(2,0.25) circle (0.25) {}
(3,0.25) circle (0.25) {}
(4,0.25) circle (0.25) {}
(5,0.25) circle (0.25) {}
(6,0.25) circle (0.25) {}
(7,0.25) circle (0.25) {}
(8,0.25) circle (0.25) {}
(9,0.25) circle (0.25) {}
(10,0.25) circle (0.25) {}
(11,0.25) circle (0.25) {}
(12,0.25) circle (0.25) {};

\draw[thick] (3,0.75) circle (0.25) {} 
(3,1.25) circle (0.25) {} 
(3,1.75) circle (0.25) {} 
(3,2.25) circle (0.25) {} 
(3,2.75) circle (0.25) {} 
(3,3.25) circle (0.25) {} 
(3,3.75) circle (0.25) {} 
(3,4.25) circle (0.25) {} 
(3,4.75) circle (0.25) {} 
(3,5.25) circle (0.25) {} 
(3,5.75) circle (0.25) {} ;

\filldraw[thick] (3,6.25) circle (0.25) {} 
(3,6.75) circle (0.25) {} 
(3,7.25) circle (0.25) {} 
(3,7.75) circle (0.25) {} 
(3,8.25) circle (0.25) {} 
(3,8.75) circle (0.25) {} 
(3,9.25) circle (0.25) {} ;

\draw[thick, opacity=0.2] (8,0.75) circle (0.25) {} 
(8,1.25) circle (0.25) {} 
(8,1.75) circle (0.25) {} 
(8,2.25) circle (0.25) {} 
(8,2.75) circle (0.25) {} 
(8,3.25) circle (0.25) {} 
(8,3.75) circle (0.25) {} ;

\draw[thick] (10,0.75) circle (0.25) {} 
(10,1.25) circle (0.25) {} 
(10,1.75) circle (0.25) {} 
(10,2.25) circle (0.25) {} ;

\draw[thick] (0,0) -- (13,0);

\draw[thick,<-] (3.75,4) to [out=30,in=150] (7.25,4);
\end{tikzpicture}    
    };
\node[inner sep=0pt, anchor= south] (configsplit) at (8,0)
    {
	\begin{tikzpicture}[scale = 0.3]
\draw (3,-.20) node[anchor=north]  {\large{$i$}}
(8,-.20) node[anchor=north]  {\large{$j$}}
;

\draw[thick, magenta] (1,0.25) circle (0.25) {} 
(2,0.25) circle (0.25) {}
(3,0.25) circle (0.25) {}
(4,0.25) circle (0.25) {}
(5,0.25) circle (0.25) {}
(6,0.25) circle (0.25) {}
(7,0.25) circle (0.25) {}
(8,0.25) circle (0.25) {}
(9,0.25) circle (0.25) {}
(10,0.25) circle (0.25) {}
(11,0.25) circle (0.25) {}
(12,0.25) circle (0.25) {};

\draw[thick] (3,0.75) circle (0.25) {} 
(3,1.25) circle (0.25) {} 
(3,1.75) circle (0.25) {} 
(3,2.25) circle (0.25) {} 
(3,2.75) circle (0.25) {} 
(3,3.25) circle (0.25) {} 
(3,3.75) circle (0.25) {} ;

\draw[thick, opacity=0.2]
(3,4.25) circle (0.25) {} 
(3,4.75) circle (0.25) {} 
(3,5.25) circle (0.25) {} 
(3,5.75) circle (0.25) {} ;

\filldraw[thick]
(6,0.75) circle (0.25) {} 
(6,1.25) circle (0.25) {} 
(6,1.75) circle (0.25) {} 
(6,2.25) circle (0.25) {} ;

\draw[thick] (8,0.75) circle (0.25) {} 
(8,1.25) circle (0.25) {} 
(8,1.75) circle (0.25) {} 
(8,2.25) circle (0.25) {} 
(8,2.75) circle (0.25) {} 
(8,3.25) circle (0.25) {} 
(8,3.75) circle (0.25) {} ;

\draw[thick] (10,0.75) circle (0.25) {} 
(10,1.25) circle (0.25) {} 
(10,1.75) circle (0.25) {} 
(10,2.25) circle (0.25) {} ;

\draw[thick] (0,0) -- (13,0);
\draw[thick,->] (3.75,2.75) to [out=30,in=150] (5.5,2.75);
\end{tikzpicture}
    };
\draw[->,thick] (config.south east) -- (7.5,2.5)
    node[midway,fill=white] {$S_{i}^{1/3}$};
\draw[->,thick] (config.south west) -- (0.5,2.5)
    node[midway,fill=white] {$M_{ij}$};
\end{tikzpicture}
  \label{fig:3a}
\end{subfigure}\hfil 
\caption{For example, consider the case $m=1$ where $\lim_{L\to\infty}w_L(0)=w(0)=0$ but $w(1)>0$. Configuration $\eta\sim \pi_{L,N}$ will not include any empty sites with high probability. 
Thus, when splitting we instead impose to split onto sites with single-occupation. Similarly, we leave a single particle behind when merging occupations of two sites. Here, $S_i$ and $M_{ij}$ denote the corresponding split- and merge-operations on configurations; we ignored the subsequent order-operation $\widehat{\cdot}$ for the sake of a clearer picture. }
\label{fig:31}
\end{figure}

Recall that under $\mu_{L,N}$ partitions consist of (at most) $L$ blocks. Since we will lift the split-merge dynamics to the space of particle configurations $\Omega_{L,N}$,
the resulting partitions should not exceed length $L$ either. However, a priori it is not clear that there is always an empty site to split on and the split-merge dynamics have to be slightly adapted to achieve this:
recall from \eqref{eq_nu_rho_equal_w} that without loss of generality we can assume that $\nu_\rho[\eta_x=n]=w(n)$.
Thus, Assumption \ref{aspt_A2_equiv_ensemble} guarantees that there exists $m\in \mathbb{N}$ such that $w(m)=\lim_{L\to\infty} w_L(m) > 0$. This implies that under $\mu_{L,N}$ a positive fraction of blocks has size $m/N$, see Lemma \ref{lem_limit_fraction_of_j_occ_sites} below. Hence, instead of leaving empty sites behind when merging, and splitting onto empty sites only, we can impose to leave a fraction of $m/N$ behind when merging and only split onto blocks of size $m/N$, cf. Figure~\ref{fig:31}. 
This will not affect the statistics on a macroscopic scale, where microscopic blocks are indistinguishable; which is why henceforward we assume $m=0$ for notational convenience. 

\begin{lemma}\label{lem_limit_fraction_of_j_occ_sites}
Let $\#_0(\eta)$ be the number of sites with zero occupation in the configuration $\eta$.
We have that 
\begin{equation*}
\pi_{L,N}
\left( 
\left( 
\frac{\#_0(\eta)}{L}
-
w(0)
\right)^2
\right)\to 0\, , \quad\text{as }\tdlim\,.
\end{equation*}
In particular, we have a weak law of large numbers and for every $\varepsilon>0$
\begin{align*}
\pi_{L,N}
\left[ 
\left|
\frac{\#_0(\eta)}{L}
-
w(0)
\right| > \varepsilon
\right] \to 0\, , \quad\text{as }\tdlim\,.
\end{align*}
\end{lemma}

\begin{proof}
By direct calculation of the second moment, it suffices to show that (using the product structure of $\pi_{L,N}$)
$$
\pi_{L,N}[\eta_x= 0, \eta_y=0]= \pi_{L-1,N}[\eta_x= 0]\,
\pi_{L,N}[ \eta_y=0]
$$ 
converges to $w(0)^2$ for all $1\leq x,y\leq L$, $x\neq y$, which holds due to \ref{aspt_A2_equiv_ensemble}.
The second statement follows immediately by Chebyshev's inequality.
\end{proof}

The operator $\mathcal{G}_{\theta}^{(N,\varepsilon)}$ approximates the corresponding continuous process acting on blocks of size larger than $\varepsilon$, which is characterised by the generator
\begin{equation*}
\begin{split}
\mathcal{G}^{(\varepsilon)}_{\theta}f(p) &:= 
\sumtwo{i,j=1}{i\neq j}^{\infty} p_i p_j
\mathds{1}\{ p_i,p_j \geq \varepsilon \}
 \left[ f(\widehat{M}_{ij}p) -f(p) \right]\\
&\qquad +
\theta
\sum_{i=1}^\infty p_i^2
\mathds{1}\{ p_i \geq 2\varepsilon \}
 \left[  \int_{\varepsilon}^{1-\varepsilon} f(\widehat{S}_{i}^{u}p) du   -f(p) \right]
\, .
\end{split}
\end{equation*}
We summarise the corresponding convergence behaviour of the generators in the following lemma:
\begin{lemma}\label{lem_conv_generators}
We have 
\begin{equation*}
\mathcal{G}_{\theta}\sind{N,\varepsilon} \to
\mathcal{G}_{\theta}\sind{\varepsilon}\, ,
\quad 
\text{as }
N\to\infty\, ,
\qquad 
\text{and}
\qquad
\mathcal{G}_{\theta}\sind{\varepsilon}
\to
\mathcal{G}_{\theta}\, ,
\quad 
\text{as }
\varepsilon\to 0\, ,
\end{equation*}
in the strong operator topology on bounded continuous functions $C_b(\overline{\nabla})$.
\end{lemma}
\begin{proof}
Let $g\in C_b (\overline{\nabla})$. To prove the first part of the statement, we need to show that
 $\|(\mathcal{G}_{\theta}\sind{N,\varepsilon} -\mathcal{G}_{\theta}\sind{\varepsilon}) g\|_{\infty}$ vanishes in the limit $N\to\infty$. 
 We will compare the two parts of each operator separately. We start with the merge term:
\begin{align*}
&
\left|
\sum_{i\neq j}^{L}p_i p_j
\mathds{1}\{ p_i,p_j \geq \varepsilon \}
 \left[ g(\widehat{M}_{ij}p) -g(p) \right] \right.\\
&\qquad \qquad \qquad-\left.
\frac{N}{N-1}\sum_{i\neq j}^L p_i p_j 
\mathds{1}\{ p_i,p_j \geq \varepsilon \}
\left[ g(\widehat{M}_{ij}p) -g(p) \right]
\right| \\
&\qquad\qquad\leq 
\frac{2 \|g\|_{\infty}}{N-1}\sum_{i\neq j}^{L}  p_i p_j
\mathds{1}\{ p_i,p_j \geq \varepsilon \}
\leq  \frac{2 \|g\|_{\infty}}{N-1}\, .
\end{align*}
And similarly for the split term:
\begin{align*}
&\left|
\sum_{i=1}^{L} p_i^2 
\mathds{1}\{ p_i \geq 2\varepsilon \}
\left[  \int_\varepsilon^{1-\varepsilon} g(\widehat{S}_{i}^{u}p) du   -g(p) \right] \right.\\
&\qquad \qquad-\left.
\frac{1}{N-1}\sum_{i=1}^L p_i  
\mathds{1}\{ p_i \geq 2\varepsilon \}
\sum_{k=\varepsilon N}^{N(p_i -\varepsilon)} 
\left[ g(\widehat{S}_{i}^{\frac{k}{Np_i}}p)   -g(p) \right]
\right| \\
&\quad \leq 
\sum_{i=1}^{L} 
p_i^2
\left|
\int_\varepsilon^{1-\varepsilon} g(\widehat{S}_{i}^{u}p) du
-
\frac{1}{(N-1)p_i}
\sum_{k=\varepsilon N}^{N(p_i -\varepsilon)}   g(\widehat{S}_{i}^{\frac{k}{Np_i}}p) 
\right|+ \frac{\|g\|_{\infty}}{N-1}
\sum_{i=1}^L p_i(1-p_i)\\
&\quad \leq 
\max_i\, \omega(g(\widehat{S}_{i}^{\bullet}p), \tfrac{1}{Np_i})
+
\frac{\|g\|_{\infty}}{N-1}
\, ,
\end{align*}
where $\omega$ is the modulus of continuity. The last inequality holds since the sum inside the absolute value approximates the Riemann-sum which in turn converges to the corresponding integral, where the error is controlled by the modulus of continuity. 
By \cite[Lemma 1]{IoTo20}, we know that 
\begin{align*}
    \|\widehat{S}_{i}^{u}p- \widehat{S}_{i}^{v}p\|_1
    \leq 2|u-v|p_i\,,
\end{align*}
which in particular implies 
\begin{align*}
\omega(\widehat{S}_{i}^{\bullet}p, \tfrac{1}{Np_i})=
\sup_{\substack{u,v\in [0,1]\\ |u-v|\leq (Np_i)^{-1}}}
    \|\widehat{S}_{i}^{u}p- \widehat{S}_{i}^{v}p\|_2
    \leq
    \sup_{\substack{u,v\in [0,1]\\ |u-v|\leq (Np_i)^{-1}}}
    \|\widehat{S}_{i}^{u}p- \widehat{S}_{i}^{v}p\|_1^{1/2}
    \leq \frac{\sqrt{2}}{\sqrt{N}}\,.
\end{align*}
This allows us to uniformly bound the moduli of continuity above, in the sense that
\begin{align}\label{eq_mod_cont}
    \omega(g(\widehat{S}_{i}^{\bullet}p), \tfrac{1}{Np_i})
    &=
    \sup_{|u-v|\leq \frac{1}{Np_i}} 
    \big|
    g(\widehat{S}_{i}^{u}p) - g(\widehat{S}_{i}^{v}p)
    \big|
    \nonumber \\ &\leq 
    \sup_{\|p'-q'\|_2\leq 
    \omega(\widehat{S}_{i}^{\bullet}p, \frac{1}{Np_i})
    } 
    \big|
    g(p') - g(q')
    \big| 
    =
    \omega\big(g,\omega(\widehat{S}_{i}^{\bullet}p, \tfrac{1}{Np_i})\big)\nonumber\\
    &\leq \omega(g, \tfrac{\sqrt{2}}{\sqrt{N}})\,,
\end{align}
where $\omega(g, \cdot)$ denotes the modulus of continuity w.r.t. $g$ on $(\overline{\nabla}, \|\cdot\|_2)$. The r.h.s. of \eqref{eq_mod_cont} vanishes due to uniform continuity of $g$. Here we used that the topology on $\overline{\nabla}$ induced by $\|\cdot\|_2$ coincides with the product topology.

Since the sum of the first estimate and second estimate multiplied by $\theta$ bounds $\|(\mathcal{G}_{\theta}\sind{N,\varepsilon} -\mathcal{G}_{\theta}\sind{\varepsilon}) g\|_{\infty}$ from above, we take the  thermodynamic limit $N/L\to\rho$ and conclude the first part of the lemma.\\
 
The second statement requires us to show that $\|(\mathcal{G}_{\theta} -\mathcal{G}_{\theta}\sind{\varepsilon}) g\|_{\infty}
\to 0$ as $\epsilon\to 0$. 
First, note that
\begin{align*}
(\mathcal{G}_{\theta} -\mathcal{G}_{\theta}\sind{\varepsilon}) g(p)
&=
\sumtwo{i,j=1}{i\neq j}^\infty p_ip_j
 \mathds{1}\{ p_i < \varepsilon \}
 \cup \{ p_j < \varepsilon \}
  \left[ g(\widehat{M}_{ij}p) -g(p) \right]\\
  &\qquad+
  \theta
\sum_{i=1}^\infty p_i^2
 \left[  \int_{0}^{1} g(\widehat{S}_{i}^{u}p) du   -g(p) \right] \\
 &\qquad -
 \theta
\sum_{i=1}^\infty p_i^2
\mathds{1}\{ p_i \geq  2\varepsilon \}
 \left[  \int_{\varepsilon}^{1-\varepsilon} g(\widehat{S}_{i}^{u}p) du   -g(p) \right]\\
 &=
 \sumtwo{i,j=1}{i\neq j}^\infty p_ip_j
 \mathds{1}\{ p_i < \varepsilon \}
 \cup \{ p_j < \varepsilon \}
  \left[ g(\widehat{M}_{ij}p) -g(p) \right]
  \\
 &\qquad +
 \theta
  \sum_{i=1}^\infty p_i^2 \mathds{1}\{ p_i <  2\varepsilon \}
 \left[  \int_{0}^{1} g(\widehat{S}_{i}^{u}p) du   -g(p) \right] \\
 &\qquad +
 \theta
 \sum_{i=1}^\infty p_i^2
\mathds{1}\{ p_i \geq  2\varepsilon \}
 \left[  \int_{0}^{1} g(\widehat{S}_{i}^{u}p) du-\int_{\varepsilon}^{1-\varepsilon} g(\widehat{S}_{i}^{u}p) du    \right]\, .
\end{align*}
The first two terms vanish by dominated convergence and for the last term we have the estimate 
\begin{equation*}
\left| \sum_{i=1}^\infty p_i^2
\mathds{1}\{ p_i {\geq}  2\varepsilon \}
 \left[  \int_{0}^{1} g(\widehat{S}_{i}^{u}p) du-\int_{\varepsilon}^{1-\varepsilon} \!\!\!\!\! g(\widehat{S}_{i}^{u}p) du    \right] \right|
\leq 
2\varepsilon
\|g\|_{\infty}
 \sum_{i=1}^\infty p_i^2
\leq 
2\varepsilon
\|g\|_{\infty}\, .
\end{equation*}
This concludes the proof.
\end{proof}

The proof of Proposition \ref{prop_pd_limit}, requires the following key observation which states that $\mu_{L,N}$ is approximately reversible w.r.t. the dynamics corresponding to $\mathcal{G}_{\theta}\sind{N,\varepsilon}$.

\begin{lemma}\label{lem_split_merge_transf_A}
For every $f,g\in C_b(\overline{\nabla})$ and $\varepsilon\in (0,1)$ we have
\begin{equation}\label{lemeq34}
\left| \mu_{L,N}(f\, \mathcal{G}_{\theta}\sind{N,\varepsilon}  g) 
-
 \mu_{L,N}(g\, \mathcal{G}_{\theta}\sind{N,\varepsilon}  f) 
 \right| \to 0
\end{equation}
in the thermodynamic limit $N/L\to\rho\geq 0$.
\end{lemma}

We postpone the proof of Lemma \ref{lem_split_merge_transf_A} until after the one of Proposition \ref{prop_pd_limit}. Now we have everything to state the proof of this section's main finding.

\begin{proof}[Proof of Proposition \ref{prop_pd_limit}:]
Due to compactness of the space $\mathcal{M}_1(\overline{\nabla})$ w.r.t. the topology induced by weak convergence, the sequence $(\mu_{L,N})$ has weak accumulation points. 
Let $\mu$ be such an accumulation point and $(\mu_{L_j,N_j})_j$ a subsequence converging to it. Then for each $\epsilon\in (0,1)$ we have
\begin{equation*}
\begin{split}
\left| \mu ( f \mathcal{G}_{\theta} g)  - \mu ( g \mathcal{G}_{\theta} f) \right|
&\leq 
\left| \mu ( f \mathcal{G}_{\theta} g)  -\mu_{L_j,N_j} ( f \mathcal{G}_{\theta}\sind{N_j,\varepsilon} g)  \right|
\\
&\qquad+
\left| \mu_{L_j,N_j} ( f \mathcal{G}_{\theta}\sind{N_j,\varepsilon} g)  - \mu_{L_j,N_j} ( g \mathcal{G}_{\theta}\sind{N_j,\varepsilon} f) \right|
 \\
&\qquad+
\left|  \mu_{L_j,N_j} ( g \mathcal{G}_{\theta}\sind{N_j,\varepsilon} f)  - \mu ( g \mathcal{G}_{\theta} f) \right|\, .
\end{split}
\end{equation*}
The middle term on the r.h.s. vanishes due to Lemma \ref{lem_split_merge_transf_A}, whereas the first term can be estimated by
\begin{equation}\label{eq_mu_approx_operators_estimate}
\begin{split}
\left| \mu ( f \mathcal{G}_{\theta} g)  -\mu_{L_j,N_j} ( f \mathcal{G}_{\theta}\sind{N_j,\varepsilon} g)  \right|
&\leq 
\left| \mu ( f \mathcal{G}_{\theta} g)  -\mu_{L_j,N_j} ( f \mathcal{G}_{\theta}g) \right| \\
&\qquad +
\left| \mu_{L_j,N_j} ( f \mathcal{G}_{\theta}g)  -\mu_{L_j,N_j} ( f \mathcal{G}_{\theta}\sind{N_j,\varepsilon} g)  \right|\, .
\end{split}
\end{equation}
Since $(\mu_{L_j,N_j})_j$ converges in distribution to $\mu$ and $\mathcal{G}_{\theta}g$ is bounded and continuous, see \cite[Lemma 4]{MWZZ02}, the first term on the r.h.s. of \eqref{eq_mu_approx_operators_estimate} vanishes as $L$ diverges. Also, by Lemma \ref{lem_conv_generators}, the second term on the right vanishes after taking the limit $L,N\to\infty$ before $\varepsilon\to 0$, since
\begin{align*}
    \left| \mu_{L_j,N_j} ( f ( \mathcal{G}_{\theta}  -  \mathcal{G}_{\theta}\sind{N_j,\varepsilon}) g  )\right|
    \leq 
    \| f\|_{\infty} \| ( \mathcal{G}_{\theta}  -  \mathcal{G}_{\theta}\sind{N_j,\varepsilon}) g \|_{\infty}\,.
\end{align*}

The same steps hold when applied to $\left|  \mu_{L_j,N_j} ( g \mathcal{G}_{\theta}\sind{N_j,\varepsilon} f)  - \mu ( g \mathcal{G}_{\theta} f) \right|$.
Overall, this yields 
$
\mu ( f \mathcal{G}_{\theta} g)  = \mu ( g \mathcal{G}_{\theta} f)
$
which finishes the proof.
\end{proof}

It only remains to state the proof of Lemma \ref{lem_split_merge_transf_A}.

\begin{proof}[Proof of Lemma \ref{lem_split_merge_transf_A}:]
First, we note that we can write $\mu_{L,N}(f \mathcal{G}\sind{N,\varepsilon}  g) $ in terms of the canonical distribution $\pi_{L,N}$:
\begin{align}\label{eq_mu_LN_termsof_pi}
&N(N-1) \mu_{L,N}(f \mathcal{G}\sind{N,\varepsilon}  g) \\
&=
\sum_{\eta \in \Omega_{L,N}} f\left( \frac{\widehat{\eta}}{N}  \right)
\sum_{i\neq j}^L \widehat{\eta}_i \widehat{\eta}_j
\mathds{1}\{ \widehat{\eta}_i ,\widehat{\eta}_j \geq \varepsilon N \}
\left[ g\left(\widehat{M}_{ij}\left(\frac{\widehat{\eta}}{N}\right) \right) - g \left(\frac{\widehat{\eta}}{N} \right) \right] \pi_{L,N}(\eta)\nonumber
\\
&+ 
\theta
\sum_{\eta \in \Omega_{L,N}}
f\left( \frac{\widehat{\eta}}{N}\right)
\sum_{i=1}^L \widehat{\eta}_i
\mathds{1}\{ \widehat{\eta}_i \geq 2\varepsilon N  \}
\sum_{k=\varepsilon N}^{\widehat{\eta}_i-\varepsilon N}
\left[ g\left(\widehat{S}_{i}^{k/\widehat{\eta}_i}\left(\frac{\widehat{\eta}}{N}\right) \right) - g \left(\frac{\widehat{\eta}}{N} \right) \right]  \pi_{L,N}(\eta)
\, ,\nonumber
\end{align}
where we multiplied with $N(N-1)$ for convenience. \\

In the following it will be easier not to work with discrete partitions, but with corresponding particle configurations without worrying about the order of the corresponding sites. 
Hence, we require a new notation to lift split and merge operations to the space of particle configurations: for $\eta\in\Omega_{L,N}$ we define
\begin{equation*}
\begin{split}
M_{xy}\eta := \eta + \eta_y (e^{(x)} - e^{(y)})\, , \\
S_{xy}^k \eta := \eta + k (e^{(y)} -e^{(x)})\, , \\
S_x^k \eta := \eta + k (e^{(L+1)}-e^{(x)})\, ,
\end{split}
\end{equation*}
where $e^{(x)} \in\Omega_{L,1}$ denotes the configuration with a single particle at site $x$, i.e. $(e^{(x)})_z =\delta_{x,z}$. 
The operator $S_x^k$ is only necessary for the case of a full particle configuration $\eta$, i.e. $\#_0 (\eta) =0$. Then, we simply append the additional block of particles at the end of the configuration, which then has length $(L+1)$. This arbitrary but convenient choice of course does not change the projection of the dynamics on the level of partitions. 
As such we can rewrite \eqref{eq_mu_LN_termsof_pi} as 
\begin{align}\label{eq_mu_as_pi_exp}
N(N-1) &\mu_{L,N}(f \mathcal{G}\sind{N,\varepsilon}  g) 
=
\pi_{L,N}\left(
 \widehat{f}\left( \eta \right)
\sum_{x\neq y}^L \eta_x \eta_y
\mathds{1}\{ \eta_x, \eta_y \geq \varepsilon N \}
\left[ \widehat{g}\left(M_{xy} \eta \right) - \widehat{g} \left(\eta\right) \right] 
\right) \nonumber
\\
&\quad+ 
\theta\,
\pi_{L,N}\left(
\mathds{1}\{ \#_0(\eta)>0 \}
\widehat{f}\left( \eta\right)
\sum_{x=1}^L \eta_x
\mathds{1}\{ \eta_x  \geq 2\varepsilon N \}\right.\nonumber\\
&\qquad\quad\times\left.
\sum_{k=\varepsilon N}^{\eta_x-\varepsilon N}
\sum_{y=1}^L
\frac{\mathds{1}\{ \eta_y = 0 \}}{\#_0(\eta)}
\left[ \widehat{g}\left(S_{xy}^{k}\eta \right) - \widehat{g} \left(\eta \right) \right] 
\right) \\
&\quad +
\theta\,
\pi_{L,N}\left(
\mathds{1}\{ \#_0(\eta)=0 \}
\widehat{f}\left( \eta\right)
\sum_{x=1}^L \eta_x
\mathds{1}\{ \eta_x  \geq 2\varepsilon N \}\right.\nonumber\\
&\qquad\quad\times\left.
\sum_{k=\varepsilon N}^{\eta_x-\varepsilon N}
\left[ \widehat{g}\left(S_{x}^{k}\eta \right) - \widehat{g} \left(\eta \right) \right] 
\right)\,  , \nonumber
\end{align}
where $\widehat{f}:= f\circ T$ and $\widehat{g}:= g\circ T$ with the ordering map $T$ given in \eqref{eq_def_T}. 
Furthermore, we can see that the terms that depend on $\widehat{f}$ and $\widehat{g}$ only through the product $\widehat{f}(\eta) \widehat{g}(\eta)$ cancel when taking the difference between $\mu_{L,N}(f \mathcal{G}\sind{N,\varepsilon}  g) $ and $\mu_{L,N}(g \mathcal{G}\sind{N,\varepsilon}  f)$ in \eqref{lemeq34}. Additionally, the very last term in \eqref{eq_mu_as_pi_exp} is negligible since, using $(\eta_x-1)/(N-1)\leq 1$ and $\sum_{x=1}^L\eta_x/N = 1$ $\pi_{L,N}$-almost surely, we have 
\begin{equation*}
\begin{split}
\frac{1}{N(N-1)}
\pi_{L,N}&
\left(
\mathds{1}\{ \#_0(\eta)=0 \}
\sum_{x=1}^L \eta_x
(\eta_x-1)
\right) 
\leq 
\pi_{L,N}
\left[
 \#_0(\eta)=0
\right]\, ,
\end{split}
\end{equation*}
which vanishes due to Lemma \ref{lem_limit_fraction_of_j_occ_sites} because $w(0)>0$.\\

To simplify notation we introduce
\begin{equation}\label{veq}
V_{f,g}\sind{L,N,\varepsilon}:=\pi_{L,N}\left(
 \widehat{f}\left( \eta \right)
\sum_{x\neq y}^L \eta_x \eta_y 
\mathds{1}\{ \eta_x,\eta_y  \geq \varepsilon N \}
\widehat{g}\left(M_{xy} \eta \right) 
\right)
\end{equation}
and 
\begin{align}
U_{f,g}\sind{L,N,\varepsilon} &:=
\theta\, 
\pi_{L,N}\left(
\mathds{1}\{ \#_0(\eta)>0 \}
\widehat{f}\left( \eta\right)
\sum_{x=1}^L \eta_x
\mathds{1}\{ \eta_x  \geq 2\varepsilon N \}\right.\nonumber\\
&\qquad\qquad\times\left.
\sum_{k=\varepsilon N}^{\eta_x-\varepsilon N}
\sum_{y=1}^L
\frac{\mathds{1}\{ \eta_y = 0 \}}{\#_0(\eta)}
\widehat{g}\left(S_{xy}^{k}\eta \right) 
\right) \, .
\label{ueq}
\end{align}
We are then left to analyse   
\begin{align*}
 N(N-1) &(\mu_{L,N}(f \mathcal{G}\sind{N,\varepsilon}  g) 
-\mu_{L,N}(g \mathcal{G}\sind{N,\varepsilon}  f))\\
&=
 V_{f,g}\sind{L,N,\varepsilon}
+U_{f,g}\sind{L,N,\varepsilon}
- V_{g,f}\sind{L,N,\varepsilon}
-U_{g,f}\sind{L,N,\varepsilon}+o(N^2)\, .
\end{align*}
The goal is to compare $V_{f,g}\sind{L,N,\varepsilon}-U_{g,f}\sind{L,N,\varepsilon}$ and $U_{f,g}\sind{L,N,\varepsilon}- V_{g,f}\sind{L,N,\varepsilon}$, respectively, and show that these differences vanish in the limit if divided by $N(N-1)$.
In order to prove this, we perform a change of measure, since both $V\sind{L,N,\varepsilon}$ and $U\sind{L,N,\varepsilon}$ are expectations with respect to $\pi_{L,N}$.

First, we note that the restriction of the merge map 
\begin{equation*}
M_{xy}\Big|_{\{\eta_y= k\}}:
\Omega_{L,N}\cap \{\eta_y= k\}
\to
\Omega_{L,N}
\end{equation*}
is injective and therefore defines a bijection between the set $\Omega_{L,N}\cap \{\eta_y= k\}$ and its image $M_{xy}(\Omega_{L,N}\cap \{\eta_y= k\}) = \Omega_{L,N}\cap \{\eta_y= 0, k\leq \eta_x\}$ with inverse $S_{xy}^k$.
Therefore, the change of measure of $\pi_{L,N}$ and its pushforward measure under $M_{xy}$ for fixed $x,y,k$ on the set $\Omega_{L,N}\cap \{\eta_y= k\}$ is given by:
\begin{equation}\label{eq_rn_deriv_pi_A}
\begin{split}
\frac{\pi_{L,N} }{\pi_{L,N} \circ M_{xy}} [\eta]
&=
\frac{\prod_{z=1}^L w_L(\eta_z)}{\prod_{z=1}^L w_L((M_{xy}\eta)_z)}
=
\frac{w_L((M_{xy}\eta)_x-k)w_L(k)}{w_L((M_{xy}\eta)_x)w_L(0)}\, .
\end{split}
\end{equation}
This will allow us to perform a change of measure in the following sense.
Fix $x,y$ and $k$, furthermore let $h_{y,k}$ be a real valued function on $\Omega_{L,N}$ with support in $\{\eta_y=k\}$. 
By definition of $M_{xy}$, we can define $\widetilde{h}_{x,y,k} = h_{y,k} \circ M_{xy}^{-1}$ which is zero outside of $\Omega_{L,N}\cap \{\eta_y= 0, k\leq \eta_x\}$. We then have 
\begin{equation*}
\begin{split}
\pi_{L,N}(h_{y,k}(\eta) )
&=
\pi_{L,N}(\widetilde{h}_{x,y,k}(M_{xy}\eta) )
=
\int_{\{\eta_y= k\}} \widetilde{h}_{x,y,k}(M_{xy}\eta) \;\pi_{L,N}[d\eta]\\
&=
\int_{\{\xi_y= 0, k\leq \xi_x\}} \widetilde{h}_{x,y,k}(\xi) \;\pi_{L,N}\circ M_{xy}^{-1} [d\xi] \\
&=
\int_{\{\xi_y= 0,  k\leq \xi_x\}} \widetilde{h}_{x,y,k}(\xi)
\frac{w_L(\xi_x-k)w_L(k)}{w_L(\xi_x)w_L(0)}
 \;\pi_{L,N} [d\xi]\, .
\end{split}
\end{equation*}
Before we apply the change of measure, we divide and multiply by $(M_{xy} \eta)_x$, decompose over $\{\eta_y=k\}$, and interchange the order of integration: 
\begin{align*}
V_{f,g}\sind{L,N,\varepsilon}&=
\sum_{x\neq y}^L
\sum_{k=1}^N
\pi_{L,N}\left(
(M_{xy} \eta)_x
 \widehat{f}\left( \eta \right)
  \widehat{g}\left(M_{xy} \eta \right)
 \frac{\eta_x \eta_y}{(M_{xy} \eta)_x}
\mathds{1}\{ k, \eta_x  \geq \varepsilon N \} 
 \right.\\
&\qquad\times
 \mathds{1}\{ \eta_y=k \} 
  \mathds{1}\{ (M_{xy} \eta)_y=0 \} 
 \bigg)\\
 &=
 \sum_{x\neq y}^L
\sum_{k=1}^N
\pi_{L,N}\left(
(M_{xy} \eta)_x
 \widehat{f}\left( S_{xy}^k M_{xy} \eta \right)
  \widehat{g}\left(M_{xy} \eta \right)
 \frac{((M_{xy} \eta)_x-k)k}{(M_{xy} \eta)_x}\right.
\mathds{1}\{ k \geq \varepsilon N \}  
 \\
&\qquad\times
\mathds{1}\{ k \leq (M_{xy} \eta)_x-\varepsilon N \} 
 \mathds{1}\{ (S_{xy}^k M_{xy} \eta)_y =k \} 
  \mathds{1}\{ (M_{xy} \eta)_y=0 \} 
 \bigg)\, ,
\end{align*}
where in the last step we expressed $\eta$ in terms of $M_{xy}\eta$.
Recall now \eqref{eq_rn_deriv_pi_A} and note that we can recover $\eta$ from $ M_{xy}\eta$ and $k=\eta_y$. Therefore, the change of measure yields
\begin{align*}
V_{f,g}\sind{L,N,\varepsilon} &=
 \sum_{x\neq y}^L
\sum_{k=1}^N
\pi_{L,N}\left(
\mathds{1}\{ \#_0(\eta)>0 \}
\eta_x
 \widehat{f}\left( S_{xy}^k  \eta \right)
  \widehat{g}\left(\eta \right)
 \frac{( \eta_x-k)k}{ \eta_x}
 \frac{w_L(\eta_x-k)w_L(k)}{w_L(\eta_x)w_L(0)}.
\right.\\
&\qquad\times
\mathds{1}\{ k \leq  \eta_x-\varepsilon N \} 
\mathds{1}\{ k \geq \varepsilon N \}  
 \mathds{1}\{ (S_{xy}^k \eta)_y =k \} 
  \mathds{1}\{ \eta_y=0 \} 
 \bigg)\\
&=
\pi_{L,N}\left(
\mathds{1}\{ \#_0(\eta)>0 \}
  \widehat{g}\left(\eta \right)
 \sum_{x=1}^L
 \eta_x
 \mathds{1}\{ \eta_x \geq 2\varepsilon N \}  
 \sum_{k=\varepsilon N}^{\eta_x-\varepsilon N}
  \sum_{y=1}^L
 \frac{ \mathds{1}\{ \eta_y=0 \}}{\#_0(\eta)} 
\right.\\
&\qquad\times
\widehat{f}\left( S_{xy}^k  \eta \right)
\#_0(\eta)
 \frac{( \eta_x-k)k}{ \eta_x}
 \frac{w_L(\eta_x-k)w_L(k)}{w_L(\eta_x)w_L(0)}
 \bigg)\, .
\end{align*}
Comparing now $V_{f,g}\sind{L,N,\varepsilon}$ to $U_{g,f}\sind{L,N,\varepsilon}$ in \eqref{ueq}, the only discrepancy is the term
\begin{equation}\label{eq_pre_punchline}
\begin{split}
\#_0(\eta)
 \frac{( \eta_x-k)k}{ \eta_x}
 \frac{w_L(\eta_x-k)w_L(k)}{w_L(\eta_x)w_L(0)}
= 
\frac{\#_0(\eta)}{L}
\frac{k\,w_L(k)L}{w_L(0)}
\frac{(\eta_x-k)\,w_L(\eta_x-k)}{\eta_x\, w_L(\eta_x) }
 \, .
\end{split}
\end{equation}
First, note that by Assumption \ref{aspt_A3_theta}
\begin{equation}\label{eq_punchline}
k\,w_L(k)L
\frac{(\eta_x-k)w_L(\eta_x-k)}{\eta_x w_L(\eta_x) }
\to \theta\, ,\quad\mbox{as }L\to\infty\, ,
\end{equation}
uniformly in $k$ and $\eta$, since $\varepsilon N \leq k, \eta_x -k, \eta_x$. 
Therefore, we have
\begin{align*}
&\left| \mu_{L,N}(f \mathcal{G}\sind{N,\varepsilon}  g) 
-
\mu_{L,N}(g \mathcal{G}\sind{N,\varepsilon}  f) 
\right| \\
&\qquad =
\frac{1}{N(N-1)}
\left|
\left(V_{f,g}\sind{L,N,\varepsilon}
-U_{g,f}\sind{L,N,\varepsilon}\right)
+\left(U_{f,g}\sind{L,N,\varepsilon}
- V_{g,f}\sind{L,N,\varepsilon}\right)
\right| + o(1)\\
&\qquad\leq
\frac{2\|f\|_{\infty} \|g\|_{\infty} \theta }{N(N-1)}
\pi_{L,N}\left(
 \sum_{x=1}^L
 \eta_x\right.
 \left.\sum_{k=\varepsilon N}^{\eta_x-\varepsilon N}
\left| \frac{\#_0(\eta)}{L}
 \frac{1}{w_L(0)}
- 1 \right|
\right) + o(1)\\
&\qquad\leq 
2\|f\|_{\infty} \|g\|_{\infty} \theta\,
\pi_{L,N}\left(
\left| \frac{\#_0(\eta)}{L}
 \frac{1}{w_L(0)}
- 1 \right| 
\right)+ o(1)\, ,
\end{align*}
where we used  \eqref{eq_punchline} before dropping all indicator functions in the first inequality.
Lastly, by Lemma \ref{lem_limit_fraction_of_j_occ_sites} and Assumption~\ref{aspt_A1_pos_limit_exists} we know that
\begin{align*}
\lim_{\tdlim}
    \pi_{L,N}\left(
\left| \frac{\#_0(\eta)}{L}
 \frac{1}{w_L(0)}
- 1 \right| 
\right)=0\, ,
\end{align*}
which completes the proof.

\if{false}{
Under $\pi_{L,N}$ we have that 
\begin{equation*}
\frac{\#_0(\eta)}{L}\simeq  w(0)
\end{equation*}
 due to Lemma \ref{lem_limit_fraction_of_j_occ_sites}.
The remaining term is taken care of by 
assumption \ref{aspt_A3_theta}, which implies that  
\begin{equation}\label{eq_punchline}
k\,w_L(k)L
\frac{(\eta_x-k)w_L(\eta_x-k)}{\eta_x w_L(\eta_x) }
\to \theta\quad\mbox{as }L\to\infty
\end{equation}
uniformly in $k$ and $\eta$, since $\varepsilon N \leq k, \eta_x -k, \eta_x$. 
Overall, 
\begin{align*}
&\left| \mu_{L,N}(f \mathcal{G}\sind{N,\varepsilon}  g) 
-
\mu_{L,N}(g \mathcal{G}\sind{N,\varepsilon}  f) 
\right| \\
&\qquad\simeq
\frac{1}{N(N-1)}
\left|
\left(V_{f,g}\sind{L,N,\varepsilon}
-U_{g,f}\sind{L,N,\varepsilon}\right)
+\left(U_{f,g}\sind{L,N,\varepsilon}
- V_{g,f}\sind{L,N,\varepsilon}\right)
\right|\\
&\qquad\leq
\frac{2\|f\|_{\infty} \|g\|_{\infty}}{N(N-1)}
\pi_{L,N}\left(
\mathds{1}\{ \#_0(\eta)>0 \}
 \sum_{x=1}^L
 \eta_x\right.\\
 &\qquad\qquad \times\left.\sum_{k=\varepsilon N}^{\eta_x-\varepsilon N}
\left| \#_0(\eta)
 \frac{( \eta_x-k)k}{ \eta_x}
 \frac{w_L(\eta_x-k)w_L(k)}{w_L(\eta_x)w_L(0)}
-\theta \right|
\right)\\
&\qquad\leq 
2\|f\|_{\infty} \|g\|_{\infty}
\pi_{L,N}\left(
\left| \#_0(\eta)
 \frac{( \eta_x-k)k}{ \eta_x}
 \frac{w_L(\eta_x-k)w_L(k)}{w_L(\eta_x)w_L(0)}
-\theta \right|
\right) \to 0,
\end{align*}
as $N/L\to\rho$, which completes the proof.
}\fi
\end{proof}
\smallskip

\subsection{Concentration and uniqueness of the limit}


As already mentioned in the introduction,
the only invariant distribution w.r.t. $\mathcal{G}_{\theta}$ which concentrates on full partitions $\nabla_{[0,\alpha]}$, for some $\alpha\in [0,1]$, is the Poisson-Dirichlet distribution PD$_{[0,\alpha]}(\theta)$. 
Hence, in order to prove that weak  accumulation points of measures $(\mu_{L,N})_{L,N}$ coincide, it is enough to show that each such weak limit $\mu$ from Proposition \ref{prop_pd_limit} satisfies
\begin{equation*}
\var_{\mu}\left(\sum_{i=1}^{\infty}p_i\right) = 0
\end{equation*}
and there exists $\alpha_\mu \in[0,1]$ such that $\mu\left(\sum_{i=1}^{\infty}p_i\right) = \alpha_\mu$ for every limit $\mu$.


In the following we consider not only the $\ell^1$-norm but also the $\ell^2$-norm as a function on $\overline{\nabla}$ and write
\begin{equation*}
\|p\|_1 = \sum_{i=1}^{\infty} p_i
\quad 
\text{ and }
\quad 
\|p\|_2^2 = \sum_{i=1}^{\infty} |p_i|^2
\, ,
\end{equation*}
respectively. Note that the latter has the advantage, in contrast to $\|\cdot\|_1$, of being a continuous function on $\overline{\nabla}$ with respect to the product topology. \\

First, we recall the first part of Lemma 5 in \cite{MWZZ02} where Mayer-Wolf et al. proved the following relationship between $\mu(\|p\|_1^2)$ and $\mu(\|p\|_2^2)$, which is a natural consequence of the balance between expected split and merge rates for the stationary distribution $\mu$.

\begin{lemma}[Mayer-Wolf et al.]\label{lem_Zerner}
Let $\theta\in (0,1]$ and let $\mu\in \mathcal{M}_1(\overline{\nabla})$ be invariant w.r.t. $\mathcal{G}_{\theta}$, then 
\begin{equation*}
\mu(\|p\|_1^2)= (1+\theta)\mu(\|p\|_2^2)\, .
\end{equation*}
\end{lemma}

Together with its counterpart in the following lemma, this immediately proves that each subsequential limit of $(\mu_{L,N})_{L,N}$ concentrates on some $\nabla_{[0,\alpha_\mu ]}$.

\begin{lemma}\label{lem_sb_mean_sq}
Let $\theta\in (0,1]$ and let $\mu \in \mathcal{M}_1(\overline{\nabla})$ be an accumulation point of  $(\mu_{L,N})_{L,N}$.
Then 
\begin{equation*}
\mu (\|p\|_1)^2
 =
(1+\theta) \mu (\|p\|_2^2) 
\, .
\end{equation*}
\end{lemma}

\begin{proof}
Since the statement is trivial for $\mu = \delta_{(0,0,\ldots)}$, we assume without loss of generality $\mu \neq \delta_{(0,0,\ldots)}$.
Let $(\mu_{L_j,N_j})_{j\in \mathbb{N}}$ be a sequence converging weakly to $\mu$.
For such $\mu$ we 
first note that
\begin{align}\label{eq_2nd_mmt_part2part}
    \mu (\|p\|_2^2)  &= 
\lim_{j\to\infty} \mu_{L_j,N_j}(\|p\|_2^2) 
= \lim_{j\to\infty} \frac{L_j}{N_j^2} \pi_{L_j,N_j}(\eta_1^2)\,,
\end{align}
since $\mu_{L_j,N_j}$ converges weakly to $\mu$.
Rewriting the r.h.s. as size-biased expectation, cf. \eqref{sbform}, yields in particular
\begin{align}\label{eq_ident_two_norm}
\mu (\|p\|_2^2)  &= 
 \lim_{j\to\infty} \frac{1}{N_j} \pi_{L_j,N_j}(\widetilde{\eta}_1) 
= \lim_{j\to\infty} \mu_{L_j,N_j} (\widetilde{p}_1)
= \lim_{j\to\infty} \sigma(\mu_{L_j,N_j}) (q_1)\, .
\end{align}
In the second to last step we used that the size-biased distribution is invariant under reordering of the configuration.
Furthermore, weak convergence of $\mu_{L_j,N_j}$ implies weak convergence of the size biased distributions by Lemma \ref{lem_size_biased_conv}, i.e. 
\begin{equation*}
\sigma(\mu_{L_j,N_j}) \to
\sigma(\mu)
\end{equation*}
in the topology induced by weak convergence. Recall that under $\sigma(\mu)$, the first component in $q$ is zero with probability 
\begin{equation}\label{eq_lim_size_biased_eqal_0}
\sigma(\mu)[q_1=0] = \mu(1-\|p\|_1) \,  .
\end{equation}


Since the projection map on the first component is continuous w.r.t. the product topology, we have $\sigma(\mu_{L,N})(q_1) \to \sigma(\mu)(q_1)$ which yields with \eqref{eq_ident_two_norm} and \eqref{eq_sb_fm}
\begin{equation}\label{intermediate}
\mu (\|p\|_2^2)  = 
\sigma( \mu) (q_1)
=
\sigma( \mu) (q_1\; | \; q_1>0) \sigma( \mu)[ q_1>0 ]
= \sigma'(\mu)(q_1) \mu(\|p\|_1)\, ,
\end{equation}
where $\sigma'(\mu)[q_1 \in \cdot ]$ is the size-biased distribution conditioned on positive components defined. 
Now, using \eqref{eq_def_pos_sb} and disintegration of measures we write 
\begin{align*}
\sigma'(\mu)(q_1)
&=
\int_{\overline{\nabla}} \|p\|_1 \int_{\nabla}
q_1 \, \sigma_{p/\|p\|_1}[dq]\,
\mu[dp]\\
&=
\int_0^1 \bar{\alpha}
\int_{\overline{\nabla}}  \int_{\nabla}
q_1 \, \sigma_{p/\|p\|_1}[dq]\,
\mu^{(\bar{\alpha})}[dp]\, \mu[\|p\|_1\in d\bar{\alpha} ]\, ,
\end{align*}
where $\mu^{(\bar{\alpha} )} =\mu [ \,\cdot\, |\| p\|_1 =\bar{\alpha}]$ denotes the measure $\mu$ conditioned on $\nabla_{[0,\bar{\alpha} ]}$, which is well defined almost surely w.r.t.\ $\mu[\|p\|_1\in d\bar{\alpha} ]$.
Since $\mathcal{G}_\theta$ \eqref{eq_discrete_gen_split_merge} conserves $\| p\|_1$, $\mu^{(\bar{\alpha} )}$ is also invariant for $\mathcal{G}_\theta$ and thus equal to PD$_{[0,\bar{\alpha} ]}(\theta )$ by Proposition \ref{prop_PD_unique_inv}. 
Hence, $\sigma(\mu^{(\bar{\alpha})})[dq]=\int_{\overline{\nabla}}\sigma_{p/\|p\|_1}[dq]\,
\mu^{(\bar{\alpha})}[dp]$, as in \eqref{eq_def_sigma_mu}, is the GEM distribution and therefore
\begin{equation*}
\int_{\overline{\nabla}}  \int_{\nabla}
q_1 \, \sigma_{p/\|p\|_1}[dq]\,
\mu^{(\bar{\alpha})}[dp] = \frac{1}{1+\theta}\, , \quad \text{ for a.e. }\bar{\alpha}\in (0,1]\, .
\end{equation*}
Thus,
\begin{equation*}
\sigma '( \mu) (q_1)
=
\frac{1}{1+\theta}
\int_0^1
\bar{\alpha}  \;\mu[ \|p\|_1\in d\bar{\alpha}]
=\frac{1}{1+\theta} \mu (\|p\|_1)
\end{equation*}
which, with \eqref{intermediate}, yields with 
\begin{equation*}
\mu (\|p\|_2^2) 
=
\frac{1}{1+\theta}  \mu (\|p\|_1)^2 \, ,
\end{equation*}
concluding the proof.
\end{proof}

\begin{remark}
It is interesting to note that Lemma \ref{lem_sb_mean_sq} holds in greater generality, because all $k$-norms are continuous whenever $k>1$. We have
\begin{align*}
    \alpha_{\mu}^k
    =
   \mu(\|p\|_1)^k
    =
    \frac{1}{(k-1)!}\left[\prod_{j=1}^{k-1}(j+\theta)\right]  \mu(\|p\|_k^k)
    \qquad\text{for all } k\geq 1\, .
\end{align*}
This can be proven in exactly the same way as Lemma \ref{lem_sb_mean_sq}, with the only difference that we replace $\mu(\|p\|_2^2)$ with $\mu(\|p\|_k^k)$ when $k>1$. The case $k=1$ is trivial.
\end{remark}

Lemmas \ref{lem_Zerner} and \ref{lem_sb_mean_sq} imply $\mu (\|p\|_1^2)=\mu (\|p\|_1)^2$, which is equivalent to
\begin{equation}\label{eq_zero_variance}
\var_{\mu}(\|p\|_1) = \mu (\|p\|_1^2)- \mu (\|p\|_1)^2
= 0\, .
\end{equation}
Therefore, accumulation points $\mu$ concentrate on $\nabla_{[0,\alpha_{\mu}]}$ with $\alpha_{\mu}:= \mu(\|p\|_1)$. To conclude the proof of Theorem \ref{THEO_MAIN}, it is left to show that $\alpha_{\mu}$ is independent of the choice of $\mu$ if assumption \ref{aspt_A4_concentration} is satisfied.

\begin{proof}[Proof of Theorem \ref{THEO_MAIN}]
Let $\theta\in (0,1]$.
By Proposition \ref{prop_pd_limit}, we know that under Assumptions \sparseapmts all weak accumulation points of $(\mu_{L,N})_{L,N}$ are reversible w.r.t. $\mathcal{G}_{\theta}$. 
Furthermore, every such limit $\mu$ of a subsequence $(\mu_{L_j,N_j})_{j\geq 1}$ concentrates on $\nabla_{[0,\alpha_{\mu}]}$ for some $\alpha_{\mu}\in [0,1]$ as follows from \eqref{eq_zero_variance}. 
Thus, Proposition \ref{prop_PD_unique_inv} implies that $\mu$ must be the Poisson-Dirichlet distribution on $[0,\alpha_{\mu}]$ with parameter $\theta$, which proves the first statement.

The only control we have on $\alpha_{\mu}$ from \sparseapmts is that
\begin{equation*}
\alpha_\mu \leq 1-\frac{\nu_{\rho}(\eta_1)}{\rho}\, , \quad\text{where}\quad\nu_{\rho}(\eta_1)  =\sum_{n=0}^{\infty}n w(n)
\end{equation*}
is the expected density in the bulk. 
However, Assumption \ref{aspt_A4_concentration} together with  Lemma \ref{lem_sb_mean_sq} implies uniqueness of the limit, since for every accumulation point $\mu$ we have
\begin{equation}\label{eq_subseq_mean_one_norm}
\mu(\|p\|_1) =  \sqrt{(1+\theta)\mu(\|p\|_2^2)}
= \sqrt{\frac{1+\theta}{\rho}\lim_{\tdlim}\frac{\pi_{L,N}(\eta_x^2)}{N} }
= \alpha\, ,
\end{equation}
where we used the identity \eqref{eq_2nd_mmt_part2part} in the second equality.
This implies the second statement and concludes the proof of the theorem. 
\end{proof}

\medskip

\section{Proof of Theorem \ref{THEO_SPECIALISED_MAIN}}\label{sec_proof_specific_result}
In this section we prove Theorem \ref{THEO_SPECIALISED_MAIN}, a specialized version of Theorem \ref{THEO_MAIN} which is better suited for application to condensation in particle systems.
In contrast to assumptions \ref{aspt_A1_pos_limit_exists} -\ref{aspt_A3_theta}, we require uniform convergence of the weights as well as stronger control on the weights in sub-$L$ scales in \ref{aspt_B3_theta}. However, this allows us to drop assumption \ref{aspt_A4_concentration} that was needed to guarantee the concentration of the macroscopic phase in Theorem \ref{THEO_MAIN}. 
 Thanks to the stronger assumptions, we can explicitly calculate the form of limiting single-site marginals, which will imply the equivalance of ensembles \ref{aspt_A2_equiv_ensemble} and  the condensation transition. 

In Appendix \ref{sec:equivensembles} we prove the equivalence of ensembles and deduce that the system defined by weights $(w_L)_L$ satisfying \ref{aspt_B1_uniform_conv} (and a growth condition on the weights which is weaker than \ref{aspt_B3_theta}) exhibits condensation for $\rho$ large enough. We summarize this result here, a more general and detailed version is given in Proposition \ref{prop_single_site_marg}.

\begin{proposition}[Equivalence of ensembles]\label{equivalence1}
Consider weights $(w_L)_L$ and $w$ satisfying \ref{aspt_B1_uniform_conv}, \ref{aspt_B3_theta} with corresponding canonical measures $(\pi_{L,N})_{L,N}$ as defined in \eqref{wpi}. 
Then the system exhibits a condensation transition in the thermodynamic limit $\tdlim\geq 0$ (cf. Definition \ref{def_cond}) with critical density $\rho_c :=\sum_{n=1}^\infty nw(n)$. 
More precisely, we have convergence of single-site marginals
$\pi_{L,N} [\eta_x \in\cdot ]\to\nu_\rho$ such that
\begin{equation}
    \nu_\rho (\eta_x )=\left\{\begin{array}{cl}
        \rho & \text{if }\rho <\rho_c\, , \\
        \rho_c & \text{if }\rho\geq \rho_c\, .
    \end{array}   \right.
\end{equation}
and $\nu_\rho =w$ for $\rho\geq \rho_c$.
\end{proposition}
This establishes existence of the condensed phase for $\rho >\rho_c$, and the following result guarantees that it agrees with the macroscopic phase, and there is no mass on intermediate scales.

\begin{proposition}\label{prop_concentration}
Consider weights $(w_L)_L$ and $w$ satisfying \ref{aspt_B1_uniform_conv} and \ref{aspt_B3_theta} with corresponding canonical measures $(\pi_{L,N})_{L,N}$.  Then, for every $\rho>\rho_c$
\begin{align*}
    \lim_{\varepsilon\to 0}\lim_{\tdlim}\pi_{L,N}[\widetilde{\eta}_1 > \varepsilon N] = 1-\frac{\rho_c}{\rho}\, .
\end{align*}
\end{proposition}

\begin{proof}
We fix a density $\rho>\rho_c$.
Now, by definition of the first size-biased marginal, cf. \eqref{sbform}, we have for every $\varepsilon\in (0,1]$
\begin{align*}
    \pi_{L,N}[J<\widetilde{\eta}_1\leq \varepsilon N]
    =
    \frac{L}{N}\sum_{n=J+1}^{\varepsilon N} n\, w_L(n) \frac{Z_{L-1,N-n}}{Z_{L-1,N}}\, ,
\end{align*}
for $N$ large enough. By Assumption \ref{aspt_B3_theta}, there exists a sufficiently large $J=J(\varepsilon)$  such that
\begin{align*}
    \lim_{L\to \infty}\sup_{n\geq J}|n\, w_L(n) L -\theta | < \varepsilon\,.
\end{align*}
On the other hand, applying Lemma \ref{lem_pf_conv_1} below yields
\begin{align*}
\limsup_{\tdlim}
   \sup_{J<n\leq \varepsilon N}\frac{Z_{L-1,N-n}}{Z_{L-1,N}} \leq \Big(1- \tfrac{\varepsilon}{1-\rho_c/\rho}\Big)^{-1}\, .
\end{align*}
Altogether, 
\begin{align*}
\lim_{\tdlim}
    \pi_{L,N}[J<\widetilde{\eta}_1\leq \varepsilon N]
    \leq \varepsilon\, (\theta+\varepsilon) \Big(1- \tfrac{\varepsilon}{1-\rho_c/\rho}\Big)^{-1}\, ,
\end{align*}
which vanishes as we take the limit $\varepsilon\to 0$.
This completes the proof, since
$$
\lim_{\tdlim}
\pi_{L,N}[\widetilde{\eta}_1\leq J]= 
\lim_{\tdlim}
\frac{L}{N} \sum_{n=0}^J n\, \pi_{L,N} [\eta_1=n ]
=  \frac{\rho_c}{\rho}\, ,
$$
 which is a direct implication of Proposition \ref{equivalence1}.
\end{proof}

We now prove the key estimate used in Proposition \ref{prop_concentration}, which guarantees that the ratio of partition functions $Z_{L-1,N-n}/Z_{L,N}$ does not blow up for $n= o(N)$.

\begin{lemma}\label{lem_pf_conv_1}
Assume that \ref{aspt_B1_uniform_conv} and \ref{aspt_B3_theta} are both satisfied. Then we have
\begin{align*}
    \limsup_{N/L\to\rho}\frac{Z_{L-1,(1-\kappa)N}}{Z_{L,N}}
    \leq 
    \Big(1- \tfrac{\overline{\kappa}}{1-\rho_c/\rho}\Big)^{-1}
    \, ,
\end{align*} 
for every $\kappa=\kappa(L)=O(1)$ such that $\kappa\leq 1-\tfrac{\rho_c}{\rho}$, where $\overline{\kappa}:= \limsup_{L\to\infty}\kappa(L)$.
\end{lemma}

\begin{proof}
Fix $\e\in (0,1)$, by \ref{aspt_B3_theta} and Remark \ref{def_A} there exists $J \in \mathbb{N}_0$ such that
\begin{align}
\label{eq:w_ratio_bnd}
    \frac{m\,w_L(m)}{n\, w_L(n)} \geq 1 - \e  \quad \textrm{for all } m,n \geq J\,,
\end{align}
for all $L$ sufficiently large (depending on $J$ and $\varepsilon$).
Let $K = \lceil \kappa N\rceil$, by definition of the canonical measures we have 
\begin{align*}
\pi_{L,N}(\eta_1\,;\, \eta_1 > J) &\geq \pi_{L,N}(\eta_1\,;\, \eta_1 > J +K)= \sum_{n>J+K}^N n\,w_L(n)\frac{Z_{L-1,N-n}}{Z_{L,N}} \nonumber\\
&= \sum_{n>J}^{N-K} (n+K)\,w_L(n+K)\frac{Z_{L-1,N-K-n}}{Z_{L,N}}\nonumber\\
&\geq \big(1 - \e \big) \sum_{n>J}^{N-K} n\,w_L(n)\frac{Z_{L-1,N-K-n}}{Z_{L,N}}\,,
\end{align*}
where the last inequality follows from \eqref{eq:w_ratio_bnd}.
The sums above are all non-empty because $N-K > J$ for $L,N$ sufficiently large since $\kappa <1$.
Multiplying and dividing by $Z_{L,N-K}$ we have
\begin{align}
\label{eq:Z_ratio_bnd}
\pi_{L,N}(\eta_1\,;\, \eta_1 > J) &\geq \big(1 - \e \big)  \frac{Z_{L,N-K}}{Z_{L,N}}\pi_{L,N-K}(\eta_1\,;\, \eta_1 > J)\,.
\end{align}
By equivalence of ensembles (Proposition \ref{equivalence1}), in the thermodynamic limit $N/L\to\rho \geq\rho_c$ the single site marginals of $\pi_{L,N}$ converge weakly to $\nu_\rho$ which has mean $\rho_c$, so
\[
\pi_{L,N}(\eta_1\,;\, \eta_1 > J) = \frac{N}{L} - \pi_{L,N}(\eta_1\,;\, \eta_1 \leq J) \to \rho - \rho_J\, , \quad \textrm{as } N/L \to \rho\,,
\]
where $\rho_J$ is given by $\nu_\rho (\eta_1\,;\,\eta_1 \leq J)$. 
Furthermore, by dominated convergence, $\rho_J \to \rho_c$ as $J\to \infty$.
Now taking the thermodynamic limit in \eqref{eq:Z_ratio_bnd}, followed by the limit $J \to \infty$, we have
\begin{align*}
    \limsup_{N/L\to\rho}\frac{Z_{L,(1-\kappa)N}}{Z_{L,N}}
    \leq 
    \frac{\rho-\rho_c}{(1-\overline{\kappa})\rho-\rho_c} \frac{1}{1-\varepsilon}\,.
\end{align*}
From the equivalence of ensembles in Proposition \ref{equivalence1}, and the identity $\pi_{L,N}[\eta_1=0]=w_L (0)Z_{L-1,N}/Z_{L,N}$, we observe that $\lim_{N/L\to\rho} Z_{L-1,N}/Z_{L,N}=1$, which completes the proof, since we can choose $\varepsilon$ arbitrarily small after taking $J\to\infty$.
\end{proof}

We are now ready to state the full proof of Theorem \ref{THEO_SPECIALISED_MAIN} which follows by putting together the statements of Theorem \ref{THEO_MAIN}, Proposition \ref{equivalence1} and Proposition \ref{prop_concentration}.

\begin{proof}[Proof of Theorem \ref{THEO_SPECIALISED_MAIN}]
The equivalence of ensembles implies the condensation transition with critical density $\rho_c$, and with its formulation in Proposition \ref{prop_single_site_marg}, also  Assumption \ref{aspt_A2_equiv_ensemble} is satisfied.
Next, we recover assumptions \ref{aspt_A1_pos_limit_exists}, \ref{aspt_A3_theta} and \ref{aspt_A4_concentration} for a fixed choice of $\rho>\rho_c$.
Clearly, \ref{aspt_A1_pos_limit_exists} is a direct implication of \ref{aspt_B1_uniform_conv}, also \ref{aspt_A3_theta} follows immediately from \ref{aspt_B3_theta}. This already yields that the subsequential limits of the corresponding measures $\mu_{L,N}$ are Poisson-Dirichlet distributions. It is only left to show that the macroscopic phase is non-trivial and indeed agrees with the condensed phase.
Recall from \eqref{eq_lim_size_biased_eqal_0} that
\begin{align*}
    \mu(\|p\|_1)= \sigma(\mu)[q_1>0]
    = \lim_{\varepsilon\to 0}\sigma(\mu)[q_1 > \varepsilon]
\end{align*}
and, as $\varepsilon\to 0$,
\begin{align*}
    \sigma(\mu)[q_1 > \varepsilon]
    &= \lim_{\tdlim} \mu_{L,N}[\widetilde{p}_1 >\varepsilon]\\
    &= 
    \lim_{\tdlim}\pi_{L,N}[\widetilde{\eta}_1 >\varepsilon N]\to 1-\frac{\rho_c}{\rho}\, ,
\end{align*} 
which follows from Proposition \ref{prop_concentration}.
Hence, with \eqref{eq_subseq_mean_one_norm} we conclude \ref{aspt_A4_concentration} with $\alpha=\mu(\|p\|_1)=1-\tfrac{\rho_c}{\rho}$. Altogether, we verified \ref{aspt_A1_pos_limit_exists}-\ref{aspt_A4_concentration}  and so can apply Theorem \ref{THEO_MAIN}.
\end{proof}

\medskip

\section{Application to interacting particle systems and conclusion}\label{sec_application}



As mentioned in the introduction, condensation transitions occur naturally and have been studied extensively for interacting particle systems, more precisely for stochastic lattice gases which model transport phenomena and conserve the number of particles. A large class of such models with state space $\Omega_{L,N}$ has been introduced in \cite{CT85} with infinitesimal generator of the form
\begin{equation}\label{gene}
\mathcal{L}\sind{L} f(\eta) = 
\sum_{x,y =1}^L p(x,y) u(\eta_x, \eta_y) \left[
f(\eta^{xy}) - f(\eta)
\right]\, ,
\end{equation}
where $f\in C_b (\Omega_{L,N})$. Whenever $\eta_x>0$, $\eta^{xy}$ denotes the configuration $\eta -e^{(x)}+ e^{(y)}$ where one particle moved from site $x$ to site $y$. Recall that $e^{(x)}_z =\delta_{x,z}$.
The jump rate $u(n,m)\geq 0$ is a non-negative function of the occupation numbers $n$ on the departure and $m$ on the target site of a particle jump, and to avoid degeneracies we assume that $u(n,m)=0$ if and only if $n=0$. $p(x,y)$ denotes an irreducible probability kernel on $\{ 1,\ldots ,L\}$ and models the geometry of the underlying lattice. Systems have been studied, e.g.\ on regular lattices in various dimensions and with different boundary conditions, here we assume that the system is closed and conserves the total number of particles $\sum_{x=1}^L \eta_x$.

For any fixed number of particles $N\in\mathbb{N}$, the operator $\mathcal{L}\sind{L}$ defines an irreducible, continuous-time Markov process on the finite state space $\Omega_{L,N}$, which therefore has a unique invariant distribution $\pi_{L,N}$. It has been established in \cite{CT85,FGS16} that this distribution is indeed of product form and spatially homogeneous, cf. \eqref{wpi}, under the  conditions:
\begin{equation}\label{misa1}
\frac{u(n+1,m)}{u(m+1,n)}
= 
\frac{u(n+1,0)}{u(1,n)}
\frac{u(1,m)}{u(m+1,0)} \qquad \forall n,m \geq 0\, ,
\end{equation}
and at least one of the following
\begin{itemize}
    \item $p(\cdot, \cdot)$ is symmetric,
\end{itemize}
or
\begin{itemize}
    \item $p(\cdot, \cdot)$ is doubly stochastic, i.e. $\sum_y \big( p(x,y)-p(y,x)\big) =0$, and
\begin{equation}\label{misa2}
u(n,m)-u(m,n) = u(n,0) -u(m,0) \qquad \forall n,m \geq 0\, .
\end{equation}
\end{itemize}
Then the stationary weights are
\begin{equation}\label{eq_ident_inv_prod_meas}
w(n) = \prod_{k=1}^n  \frac{u(1,k-1)}{u(k,0)} \qquad \forall n\geq 0\, ,
\end{equation}
which depend only on the jump rates $u$ but not on the kernel $p$. 

For the special case of zero-range dynamics, \eqref{eq_ident_inv_prod_meas} simplifies further since $u(k,\cdot)=u(k)$, and \eqref{misa1} and \eqref{misa2} are fulfilled. In this case \eqref{eq_ident_inv_prod_meas} leads to the 
simple identification between stationary weights and rates
\begin{align}\label{eq_rates_in_terms_of_weights}
    u(n) = \frac{w(n-1)}{w(n)} \qquad \forall n\geq 1\, .
\end{align}
Due to this simple one-to-one correspondence between weights and transition rates, zero-range processes provide a generic framework of studying condensation transitions in interacting particle systems.\\

Another less restrictive simplification is to assume that $u(n,m)=u_1 (n)u_2 (m)$ is of product form, which automatically satisfies \eqref{misa1} and always leads to factorized stationary measures under symmetric dynamics, i.e. $p(\cdot, \cdot)$ is a symmetric transition kernel. The weights from \eqref{eq_ident_inv_prod_meas} now take the form
\begin{equation}\label{eq_ident_inv_prod_meas2}
w(n) = \bigg(\frac{u_1 (1)}{u_2 (0)}\bigg)^n \prod_{k=1}^n  \frac{u_2 (k-1)}{u_1 (k)} \qquad \forall n\geq 0\, .
\end{equation}
Note that due to the conservation law the exponential factor $\Big(\frac{u_1 (1)}{u_2 (0)}\Big)^n$ is usually omitted, since it cancels in the definition of $\pi_{L,N}$ \eqref{wpi}. 
One particular example is the inclusion process with rates
\begin{equation}\label{eq_indyn}
    u(n,m)=n\, (d+m)\quad\text{with }d>0\,  ,
\end{equation}
which leads to the stationary weights \eqref{eq_ip_weights}.
If we set
\begin{equation}\label{eq_dl}
    d=d(L)>0\quad\text{with}\quad dL\to\theta >0\text{ as }L\to\infty\, ,
    \end{equation}
the system exhibits a condensation transition with $\rho_c =0$ and a Poisson-Dirichlet structure with PD$_{[0,1]} (\theta )$, see  \cite[Theorem 1]{JCG19} or \eqref{eq_thm1_JCG} above. Theorem \ref{THEO_SPECIALISED_MAIN} recovers this result for $\theta\in (0,1]$. The restriction of $\theta$ is solely due to the fact that the Poisson-Dirichlet distribution is (so far) only proven to be the unique invariant distribution for the split-merge process if $\theta\in (0,1]$.\\

\begin{figure}[H]
    \centering 
\begin{subfigure}{0.5\textwidth}
  \includegraphics[width=\linewidth]{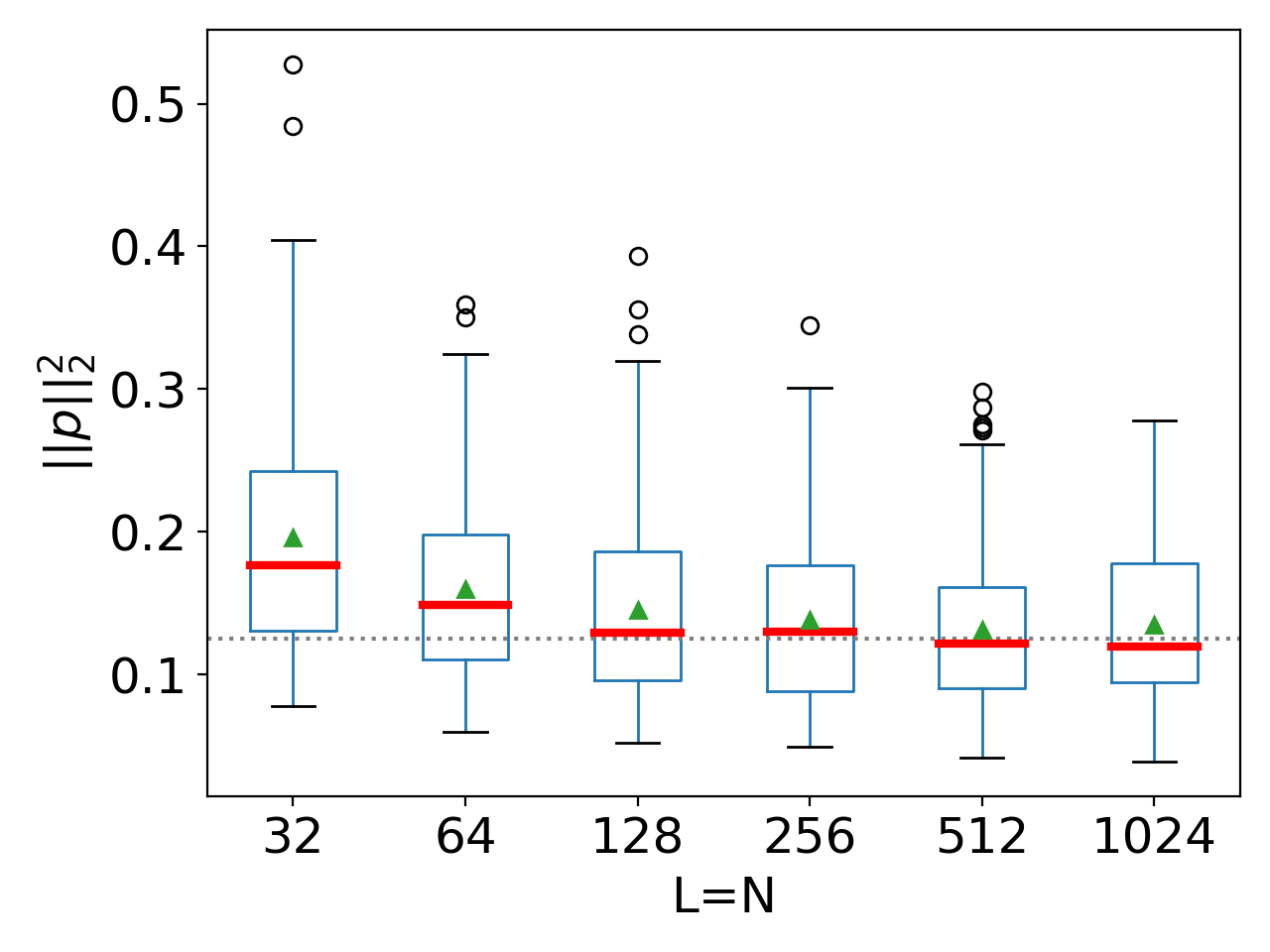}
  \label{fig:21}
\end{subfigure}\hfil 
\begin{subfigure}{0.5\textwidth}
  \includegraphics[width=\linewidth]{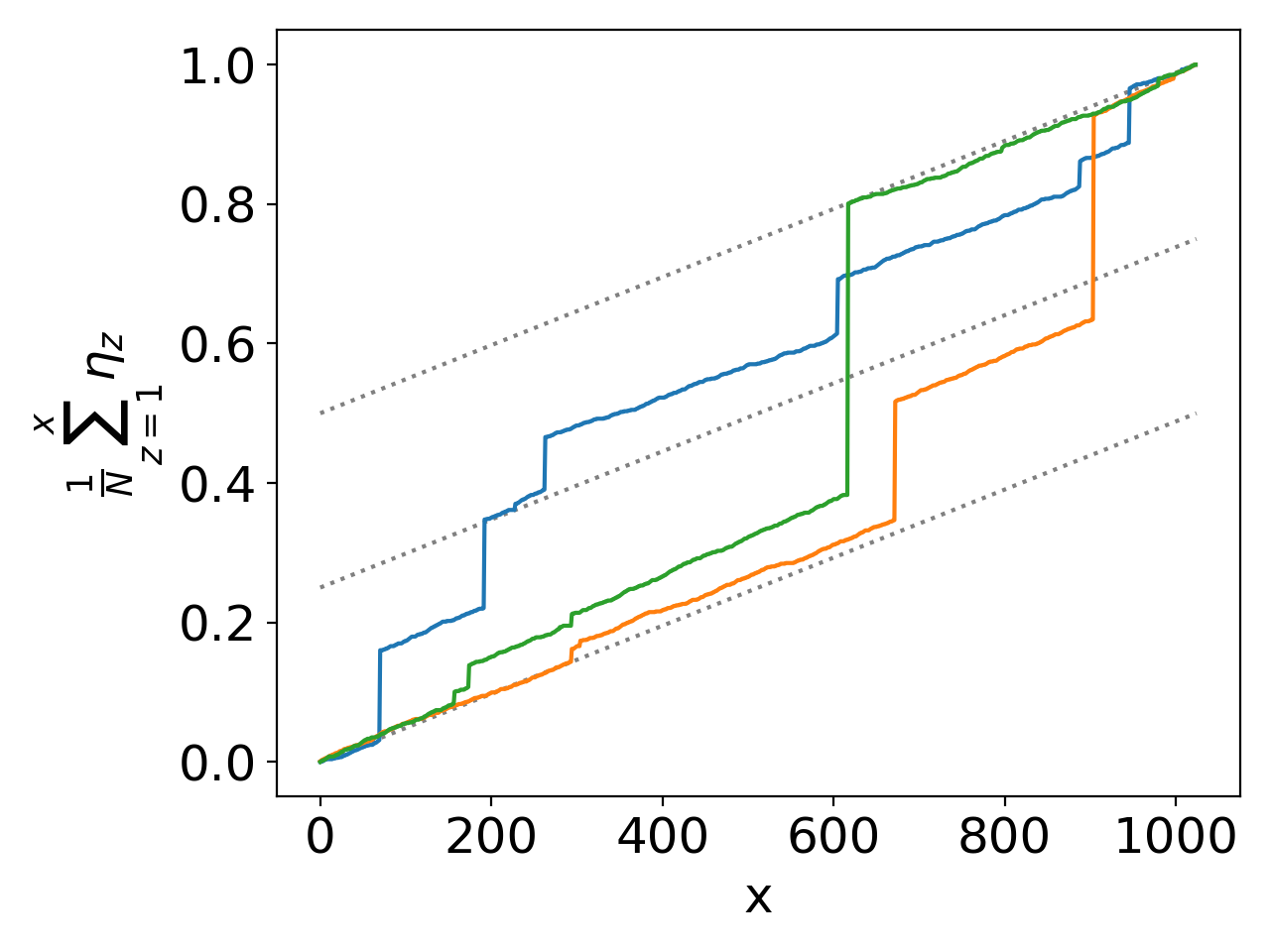}
  \label{fig:22}
\end{subfigure}\hfil 
\caption{Simulation results at stationarity for a zero-range process with rates \eqref{eq_example_rates}, parameters $\theta=A=\rho=1$ and $\rho_c =1/2$.
The left boxplot displays $\|p\|_2^2$ where $p=\tfrac{\widehat{\eta}}{N}$ was sampled $200$-times for every system size, indicating convergence to $\tfrac{1}{8}$ which agrees with $\mu_{L,N}(\|p\|_2^2)\simeq  (1+\theta)^{-1}(1-\tfrac{\rho_c}{\rho})^2$ from Lemma \ref{lem_sb_mean_sq}.
The right plot shows three samples of accumulated configurations in a system of size $N=L=1024$, where we clearly see occurance of large clusters, with a background density of $\rho_c =1/2$ indicated by dotted lines.}
\label{fig:images}
\end{figure}

Using the size-dependent parameter $d$ as above, we will provide a few instructive examples of particle systems with size-dependent jump rates $u_L$ of product form, where Theorem \ref{THEO_SPECIALISED_MAIN} applies with a non-trivial critical density $\rho_c >0$. 
We start by fixing the weights 
\begin{equation}\label{eq_example_weights}
w_L(n) =
\begin{cases}
w(n) &\quad \text{ if }n\leq A\, ,\\
\frac{d}{n} &\quad \text{ if }n> A\, ,
\end{cases}\quad\text{for some fixed }A\in\mathbb{N}\, ,
\end{equation}
where $d=d(L)$ as in \eqref{eq_dl}, $\theta\in (0,1]$ and $w$ is a probability mass function on the set $\{0,\ldots, A\}$. Let us fix for simplicity the uniform distribution with $w(n)=(1+A)^{-1}\mathds{1}\{n\leq A\}$. 
Using \eqref{eq_rates_in_terms_of_weights}, the corresponding rates of a zero-range process are given by $u_2 (m)=1$ and
\begin{equation}\label{eq_example_rates}
u_L(n) =u_{1,L} (n) =
\begin{cases}
1 &\quad \text{ if }1\leq n\leq A\, ,\\
d^{-1} \simeq L/\theta &\quad \text{ if }n= A+1\, ,\\
\frac{n}{n-1} &\quad\text{ if }n>A+1\, .
\end{cases}
\end{equation}
This underlines the mechanism that leads to condensation in such systems: Sites with occupation numbers different from $A+1$ are stable and eject particles at rates of order $1$. Sites with occupation number $A+1$ are unstable and eject a particle at diverging rate of order $L$, creating a sharp threshold between bulk sites ($\eta_x \leq A$) and cluster sites ($\eta_x >A+1$). Note that with \eqref{eq_ident_inv_prod_meas2} the same effect could be achieved by a vanishing rate of arrival onto sites with occupation number $A$.  If in the thermodynamic limit
$$
N/L\to\rho> \rho_c =\frac{A}{2}\,  ,
$$
i.e.\ the total density exceeds the expectation of the uniform bulk distribution, the system exhibits a condensation transition with mass fraction $1-\rho_c /\rho$ in the condensate. 
Heuristically, the excess mass is expelled from the bulk and accumulates in stable clusters with occupation numbers larger than $A+1$. The particular form of the rates for those clusters leads to a macroscopic phase with Poisson-Dirichlet statistics PD$_{[0,1-\rho_c /\rho]}(\theta )$, which follows from Theorem \ref{THEO_SPECIALISED_MAIN} and is illustrated in Figure~\ref{fig:images}. Clearly, the weights $(w_L)_L$ converge uniformly to $w$ and therefore satisfy both Assumption \ref{aspt_B1_uniform_conv} and \ref{aspt_B3_theta}. The same process with rates $d u_L (n)$ corresponds to asymptotically vanishing exit rates from stable sites, which is simply a time change and leads of course to the same stationary behaviour.

We can also generalize the inclusion process dynamics \eqref{eq_indyn} to stationary weights of the form \eqref{eq_example_weights}. Consider a process with rates of product form $u_L (n,m)=u_{1,L} (n) u_{2,L} (m)$ with
\begin{equation}\label{eq_indyn2}
    u_{1,L} (n)=u_{2,L} (n)=
\begin{cases}
d \simeq \theta /L&\quad \text{ if }1\leq n\leq A\, ,\\
n &\quad \text{ if }n> A\, ,
\end{cases}
\end{equation}
and $u_{2,L} (0)=d$. 
It is easy to see from \eqref{eq_ident_inv_prod_meas2} that the stationary weights for this process are given by  \eqref{eq_example_weights} and Theorem \ref{THEO_SPECIALISED_MAIN} applies. Heuristically, sites with occupation number up to $A$ eject and attract particles at a slow rate $d$, and particles on sites with higher occupation numbers become ``free'' and leave independently at the same rate and also attract other particles, as a simple generalization of the standard inclusion interaction. 
As a result, 
the dynamics in the condensed phase happen at a much higher rate than in the bulk. This separation of time scales leads to completely different dynamics than in the zero-range example above, even though both models share the same stationary distributions. We want to stress that due to the general nature of Assumptions \ref{aspt_B1_uniform_conv} and \ref{aspt_B3_theta}, Theorem \ref{THEO_SPECIALISED_MAIN} applies also to modifications of these examples and the particular form of the weights \eqref{eq_example_weights} is not important.\\



Understanding the dynamics of the condensed phase in these models is a very interesting question for future research, in particular the coarsening regime, where macroscopic clusters emerge from homogeneous initial conditions and approach stationarity by exchanging particles. Heuristically, their stationary mass partition can be understood as a balance between aggregation and fragmentation of macroscopic clusters. 
Note that these dynamics are not described by split-merge processes, which we only use as an auxiliary tool to characterize PD distributions, but are rather of a diffusive nature. A diffusive model on partitions that has stationary PD distribution has been introduced in \cite{EtKu81}, and it would be very interesting to study hydrodynamic scaling limits in this context.

As we have seen, the Poisson-Dirichlet structure in the macroscopic phase arises due to uniform stationary weights under size-biased sampling, which leads to particular rates in the zero-range process \eqref{eq_example_rates} or the generalized inclusion process \eqref{eq_indyn2} for large occupation numbers. In the context of the dynamics of interacting particle systems this is only one particular case, and it would be interesting to study the statistics of the condensed phase beyond Poisson-Dirichlet under a different scaling behaviour of the weights.


\medskip

\appendix
\section{Equivalence of ensembles}\label{sec:equivensembles}
In this section we prove that, under assumption \ref{aspt_B1_uniform_conv} and if the weights $(w_L)_L$ decay sub exponentially, i.e.
\begin{align}\label{eq_sub_exp_tails}
\frac{1}{L}\log w_L(aL)\to 0\,, \quad \forall a>0\,, \quad \text{ and } \quad
\sum_{n=0}^\infty n^2 w(n)<\infty\,,
\end{align}                                                                    
we have equivalence of ensembles and condensation in the sense of weak convergence of finite dimensional marginals. 
Note that this includes condensation in the sense of Definition \ref{def_cond}, and \eqref{eq_sub_exp_tails} is weaker than Assumption \ref{aspt_B3_theta} in Theorem \ref{THEO_SPECIALISED_MAIN}.
The result is in the same spirit as previous results on equivalence of ensembles and condensation in stochastic particle systems with stationary product measures (see for example \cite{ChGr13}).
However, as far as we know, this is the first general result in this direction for models with size-dependent weights. 

In order to state the result in more generality we first introduce some extra notation. 
For a sequence of non-negative, non-trivial weights, $(w_L(n))_{n\in\mathbb{N}_0}$, possibly depending on the system size $L$, we define a family of probability measures on $\mathbb{N}_0$ by tilting the weights by a non-negative \emph{fugacity} parameter $\phi\geq 0$:
\begin{align*}
\bar{\nu}_{\phi,L} [dn] := \frac{1}{z_L(\phi)} w_L(n) \phi^n  dn\quad\mbox{with}\quad z_L(\phi) := \sum_{n=0}^{\infty}w_L(n) \phi^n \,,
\end{align*}
which is well defined for each $\phi \in D_L := \{\phi\,:\, z_L(\phi)<\infty\}$.
The corresponding family of \emph{grand-canonical distributions} is given by the product measures 
\begin{equation}\label{eq_gc_prod_meas}
\bar{\nu}_{\phi,L}^{\otimes L } [d\eta] = \frac{1}{z_L(\phi)^L} \prod_{x=1}^L w_L(\eta_x) \phi^{\eta_x} d\eta\,,
\end{equation}
which are defined on the configuration space $\Omega_L = \mathbb{N}_0^L =\bigcup_{N=0}^\infty\Omega_{L,N}$, where the total number of particles is arbitrary. 
The expected number of particles per site (density) $R_L :D_L \to [0,\infty )$ is a strictly increasing function of $\phi$, 
with
\begin{align}
\label{eq:RL}
    R_L(\phi) := \bar{\nu}_{\phi,L} (\eta_x) = \phi\, \partial_\phi \log z_L(\phi)\, , \quad \textrm{for } \phi \in D_L\,.
\end{align}
We denote the inverse of $R_L$ by $\Phi_L$.
Furthermore, we define the variance of $\bar{\nu}_{\phi,L}$
\begin{align*}
   \sigma_L^2(\phi)=\bar{\nu}_{\phi,L}(\eta_x^2) - R_L(\phi)^2 \,,
\end{align*}
which is finite for each $\phi$ in the interior of $D_L$.
By construction, the canonical measures \eqref{wpi} on $\Omega_{L,N}$ are given by conditioning any grand-canonical measure on the total number of particles, i.e.
\begin{equation*}
\pi_{L,N} [ d\eta]  := \bar{\nu}_{\phi,L}^{\otimes L } \left[d\eta 
\;\middle|\; 
\sum_{x=1}^L  \eta_x = N \right] =  \frac{1}{Z_{L,N}} \prod_{x=1}^L w_L(\eta_x) d\eta\,,
\end{equation*}
which is independent of $\phi \in D_L$.

By Assumption \ref{aspt_B1_uniform_conv}, the weights $(w_L(n))_{n,L}$ converge uniformly in $n$ as \hbox{$L \to \infty$} to a probability measure $(w(n))_{n\in\bbN_0}$. We define the limiting grand-canonical measures by
\begin{equation*}
\bar{\nu}_{\phi} [dn] := \frac{1}{z(\phi)} w(n) \phi^n  dn\quad\mbox{with}\quad z(\phi) := \sum_{n=0}^{\infty}w(n) \phi^n\ .
\end{equation*}
Since $(w(n))_{n\in\bbN_0}$ is normalised, these measures must exist at least for each $\phi \in [0,1]$, and $\bar{\nu}_1$ corresponds to the weights $w$.
By analogy with \eqref{eq:RL}, we define the function $R(\phi) = \bar{\nu}_\phi(\eta_1)$, which is a strictly increasing function $R\colon [0,1]\to [0,\rho_c]$ with
$$
\rho_c =R(1)=\sum_{n=0}^\infty nw(n)\, ,\quad\mbox{as given in Theorem \ref{THEO_SPECIALISED_MAIN}}\ .
$$
We denote the inverse of $R$ by $\Phi$, so that the average particle density under $\bar{\nu}^{\otimes{L}}_{\Phi(\rho)}$ is $\rho$ for all $\rho\leq \rho_c$. 
Further, we denote the variance of $\bar{\nu}_\phi$ by $\sigma^2(\phi)$ which is finite for $\phi \in [0,1]$ by the second moment condition in \eqref{eq_sub_exp_tails}. Note that it may be possible that $R(\phi ) <\infty$ and $\bar\nu_\phi$ is well defined also for $\phi>1$, but such measures are not accessible as limits of $\bar\nu_{\phi ,L}$ and do not play a role in the following.

Under Assumption \ref{aspt_B1_uniform_conv} and sub-exponential weights \eqref{eq_sub_exp_tails}, it turns out that there is a condensation transition according to Definition \ref{def_cond} with critical density $\rho_c = R(1)$.
In this case, in the thermodynamic limit $\tdlim$, all finite dimensional marginals of the canonical measures converge weakly to the limiting grand-canonical measures with density $\rho$ if $\rho \leq \rho_c$, and with density $\rho_c$ if $\rho\geq\rho_c$. 
This implies that for $\rho >\rho_c$ the excess mass must condense on a vanishing volume fraction in the thermodynamic limit. 

\begin{proposition}[Equivalence of ensembles]\label{prop_single_site_marg}
Consider non-negative weights $(w_L)_L$ satisfying \ref{aspt_B1_uniform_conv} with limit $w$. Furthermore, assume that $w_L$'s have sub-exponential tails, in the sense of \eqref{eq_sub_exp_tails}.
Then for each $M \subset \mathbb{N}_0$, with $|M|=m<\infty$, denoting the marginal of $\pi_{L,N}$ on $M$ by $\pi_{L,N}^M$, we have
\begin{align*}
\pi_{L,N}^M \weakconv
\begin{cases}
\bar{\nu}_{\Phi(\rho)}^{\otimes m}&\quad \text{ if } \rho <\rho_c\,,\\
 \bar{\nu}_1^{\otimes m} &\quad \text{ if } \rho \geq \rho_c\,,
\end{cases}
\end{align*}
where $\rho_c = R(1)<\infty$.
\end{proposition}

Throughout the proof we assume further that condition \eqref{eq:bern_ass} is satisfied, i.e.
\begin{align*}
 w(0)>0\quad\mbox{and}\quad\sup_{n}[ w(n-1)\wedge w(n)]>0\,,
\end{align*}
which implies that the variance, $\s^2(\phi)$, given by $\sum n^2 w(n)\phi^n$, is positive for each $\phi \in (0,1]$. 
The special case of $w(0) =1$ is covered at the end of the proof.

\begin{figure}[H]
    \centering 
\begin{subfigure}{0.5\textwidth}
  \includegraphics[width=\linewidth]{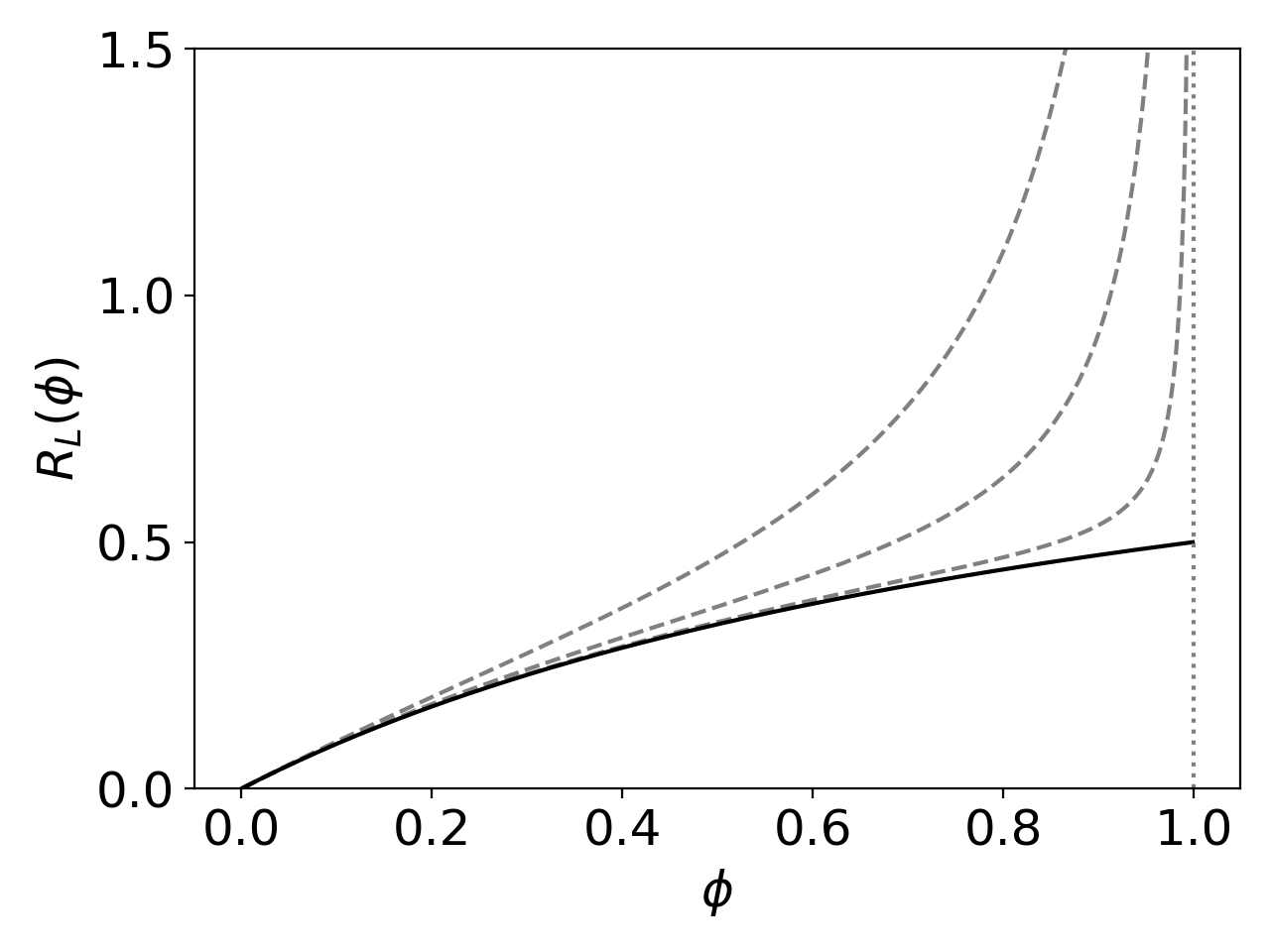}
  \label{fig:3}
\end{subfigure}\hfil 
\caption{\label{fig:typR}
The dashed lines represent $R_L(\phi)$ for systems of sizes $L\in\{4,16,128\}$ and weights \eqref{eq_example_weights} with $\theta=A=1$. The black line denotes the limit $R(\phi)=\tfrac{\phi}{1+\phi}$.
}
\label{fig:images2}
\end{figure}

We firstly observe that for each $\rho \in [0,\rho_c]$
we can construct a sequence of size-dependent fugacities such that the mean of the size-dependent grand-canonical measures converges to $\rho$ and the variance remains bounded.
The typical behaviour of $R_L(\phi)$ is illustrated in Figure \ref{fig:typR}.

\begin{lemma}
\label{lem:equiv_pre}
 Under Assumption \ref{aspt_B1_uniform_conv} and if the weights $(w_L)_{L}$ are sub-exponential in the sense of \eqref{eq_sub_exp_tails}, then for each $\varphi \in [0,1]$ there exists a sequence $(\varphi_L)_L$ in $[0,1)$, with limit point $\varphi$, such that
\begin{align}\label{eq_choice_phi}
 \|w_L - w\|_{\infty} \frac{1}{(1-\varphi_L)^3}\to 0\, , \quad \text{ as }L\to\infty\, ,
\end{align}
and
\begin{align}\label{eq_moment_limits}
    R_L(\varphi_L)\to R(\varphi)\in [0,\rho_c]\,, \quad \sigma^2_L(\varphi_L)\to \sigma^2(\varphi)>0\,,\ \text{ and }\  z_L(\varphi_L)\to z(\varphi)\,.
\end{align}
\end{lemma}
\begin{proof}
For $\varphi \in [0,1)$ we may choose $\varphi_L = \varphi$ for each $L$.
For $\varphi = 1$, we let $\varphi_L =1- \|w_L - w\|_{\infty}^{1/4}$. 
Then $\varphi_L \to 1$ as $L\to \infty$ and condition \eqref{eq_choice_phi} is satisfied since $\|w_L - w\|_{\infty} \to 0$ by assumption \ref{aspt_B1_uniform_conv}.

It is only left to prove the convergence in \eqref{eq_moment_limits}, which is equivalent to showing that
$$
\lim_{L\to\infty}\sum_{n=0}^{\infty}n^p w_L(n) \varphi_L^n = \sum_{n=0}^{\infty}n^p w(n) \varphi^n\,,
$$
for $p\in \{ 0,1,2 \}$.
Since the weights converge uniformly, we have 
\begin{align*}
    \bigg|
    \sum_{n=0}^{\infty}n^p (w_L(n)-w(n)) \varphi_L^n
    \bigg|
    \leq  \|w_L(n) -w(n)\|_{\infty}
    \sum_{n=0}^{\infty} n^p
     \varphi_L^n\, .
\end{align*}
Term-by-term differentiation of the geometric series yields for $p\in \{ 0,1,2 \}$
\begin{align*}
    \sum_{n=0}^{\infty} n^p
     \varphi_L^n = \frac{(p \, \vee 1)\varphi_L^p}{(1-\varphi_L)^{p+1}}
     + \frac{\varphi_L}{(1-\varphi_L)^{2}}\mathds{1}\{p=2\}\,.
\end{align*}
Hence, using again \eqref{eq_choice_phi}, we have
\begin{align*}
     \sum_{n=0}^{\infty}n^p w_L(n) \varphi_L^n
     \simeq 
     \sum_{n=0}^{\infty}n^p w(n) \varphi_L^n
     \to \sum_{n=0}^{\infty}n^p w(n) \varphi^n\, ,
\end{align*}
where we used dominated convergence in the last step. 
\end{proof}

We will prove Proposition \ref{prop_single_site_marg} by showing that the relative entropy between the single site marginal of $\pi_{L,N}$ and $\bar{\nu}_{\varphi_L,L}$ vanishes, where the limit density $R(\varphi)$ is equal to $\rho$ in the sub-critical case, and $\rho_c$ in the super critical case. 
Optimally, in the super-critical case, we would like to measure the relative entropy w.r.t. the limiting measure directly, however this is not possible since $\pi_{L,N}[\eta_x\in \cdot\,] \ll w$ is in general not satisfied and the relative entropy would be infinite.

The main tool we rely on in the proof of Proposition \ref{prop_single_site_marg} is a local central limit theorem (see for example \cite[Theorem 1.2]{DaMcD95}) which allows us to estimate the decay of the relative entropy. 
For completeness we include the local limit theorem here. To state it, we first introduce the Bernoulli part decomposition $q$ of a probability measure $\mathds{P}_{x,L}$ on $\Z$ as
\begin{align*}
    q(\mathds{P}_{x,L}):= \sum_{n\in\Z} (\mathds{P}_{x,L}[n]\wedge \mathds{P}_{x,L}[n+1])\, .
\end{align*}
Moreover, for a family of measures $(\mathds{P}_{x,L})_{1\leq x\leq L}$ we  define $Q_L:= \sum_{x=1}^L q(\mathds{P}_{x,L})$. 

\begin{lemma}[{\cite[Theorem 1.2]{DaMcD95}} Local central limit theorem]\label{lem_llt}
Consider a triangular array of independent integer valued random variables $\eta_{x,L}$, for $1\leq x\leq L$, and $L\in \N$, where $\eta_{x,L}$ has law $\mathds{P}_{x,L}$. Suppose there exist sequences $a_L$ and $b_L$, $L\geq 1$, such that
$b_L\to \infty$,  $\limsup_{L\to\infty} b_L^2/Q_L <\infty$
    and 
    \begin{align}\label{eq_asmpt_N01_limit}
        \frac{1}{b_L} \sum_{x=1}^L (\eta_{x,L}-a_L )\weakconv \mathcal{N}(0,1)\, .
    \end{align}
Then 
$$
\sup_{n\in \mathbb{Z}}\left|
b_L \, \mathds{P}_{L} \Big[\sum_{x=1}^L\eta_{x,L} = n \Big] -g\Big( \frac{n-a_L}{b_L}\Big)
\right|\to 0\, ,
$$
where $\mathds{P}_{L}$ denotes the product measure $\bigotimes_{ x=1 }^L\mathds{P}_{x,L}$ and $g$ the density of a standard normal.
\end{lemma}

To apply Lemma \ref{lem_llt} we consider 
independent random variables 
\begin{align}\label{eq_choice_triangular}
(\eta_{x,L})_{1\leq x\leq L}\quad \text{ with law }\quad \bar{\nu}_{\varphi_L,L}^{\otimes L}\quad\mbox{for each }L\geq1\, ,
\end{align}
where $(\varphi_L)_{L\in\mathbb{N}}$ is a sequence in $[0,1)$ satisfying \eqref{eq_choice_phi}.
To apply Lemma \ref{lem_llt} in the proof of Proposition \ref{prop_single_site_marg}, we first verify the central limit theorem \eqref{eq_asmpt_N01_limit} for the $\eta_{x,L}$'s.

\begin{lemma}\label{lem_check_lindeberg}
Consider $(\varphi_L)_L$ to be a sequence in $[0,1)$ with limit point $\varphi\in[0,1]$ satisfying \eqref{eq_choice_phi} and \eqref{eq_moment_limits}.
Furthermore, let $(\eta_{x,L})_{1\leq x\leq L, L\in \N}$ be as in \eqref{eq_choice_triangular} and define 
$$
\zeta_{x,L} := \frac{\eta_{x,L}-R_L(\varphi_L)}{\sqrt{L\, \sigma_L^2(\varphi_L)}}\, .
$$
Then $\sum_{x=1}^L\zeta_{x,L}$ converges weakly to a standard normal, as $L$ tends to infinity.
\end{lemma}

\begin{proof}
We want to apply the Lindeberg-Feller central limit theorem, see \cite[Theorem 5.12]{Ka02}: because the $\zeta_{x,L}$'s are centered and normalised, it suffices confirm that the following Lindeberg condition holds:
\[
    \textrm{For every $\varepsilon>0$, } \lim_{L\to\infty}\sum_{x=1}^L \bar{\nu}_{\varphi_L,L}(\zeta_{x,L}^2 \mathds{1}\{ |\zeta_{x,L}|>\varepsilon\})= 0\,.
\]
Since  $R_L(\varphi_L)$ and $\sigma^2_L(\varphi_L)$ converge to positive numbers, we have that for $L$ large enough 
\begin{align*}
    &\sum_{x=1}^L \bar{\nu}_{\varphi_L,L}(\zeta_{x,L}^2 \mathds{1}\{ |\zeta_{x,L}|>\varepsilon\})\\
    &= 
    \frac{1}{z_L(\varphi_L)\,\sigma^2_L(\varphi_L)}
    \sum_{n=0}^{\infty} (n-R_L(\varphi_L))^2
    w_L(n) \varphi_L^n
    \mathds{1}\Big\{ |n-R_L(\varphi_L)|>\varepsilon \sqrt{L \, \sigma^2_L(\varphi_L)}\Big\}\\
    &\leq 
    \frac{1}{z_L(\varphi_L)\,\sigma^2_L(\varphi_L)}
    \sum_{n=K_{\varepsilon,L}}^{\infty} n^2
    w_L(n) \varphi_L^n\, ,
\end{align*}
where $K_{\varepsilon,L}=\varepsilon \sqrt{L \,\sigma^2_L(\varphi_L)}+R_L(\varphi_L)$, which diverges like $\sqrt{L}$.
The denominator in the above expression converges to a positive constant.
For the numerator, we observe
\begin{align}\label{eq_conv_van_sec_mom}
     \bigg| 
    \sum_{n=K_{\varepsilon,L}}^{\infty} n^2
    (w_L(n) -w(n)) \varphi_L^n
    \bigg|
    \leq \|w_L(n) -w(n)\|_{\infty}
    \sum_{n=K_{\varepsilon,L}}^{\infty} n^2
     \varphi_L^n\to 0\, ,
\end{align}
where we used property \eqref{eq_choice_phi} of the sequence $(\varphi_L)_L$. 
Thus, using the second-moment assumption on $w$ in \eqref{eq_sub_exp_tails}, we have, for each $\epsilon$
\begin{align*}
    \sum_{n=K_{\varepsilon,L}}^{\infty} n^2
    w_L(n) \varphi_L^n
    \simeq  
    \sum_{n=K_{\varepsilon,L}}^{\infty} n^2
    w(n) \varphi_L^n \leq \sum_{n=K_{\varepsilon,L}}^{\infty} n^2
    w(n) \to 0\,,
\end{align*}
as $L\to\infty$, since $K_{\epsilon ,L} \to\infty$.
This concludes the Lindeberg condition. 
\end{proof}

\begin{proof}[Proof of Proposition \ref{prop_single_site_marg}]
We first consider the sub-critical and critical case together, fix $\rho \in (0, \rho_c]$. 
Let $(\varphi_L)_L$ be a sequence converging to $\varphi:=\Phi(\rho) \in (0,1]$ satisfying \eqref{eq_choice_phi} and \eqref{eq_moment_limits}.
We will measure the relative entropy between single-site marginals of $\pi_{L,N}[\eta_x\in\cdot\,]$ and $\bar{\nu}_{\varphi_L,L}$. 


We start with an expression for the relative entropy between $\pi_{L,N}$ and $\bar{\nu}_{\varphi_L,L}^{\otimes L}$ which is used frequently in the proof of similar equivalence of ensembles results (see for example \cite{ChGr13}), 
\begin{align*}
\begin{split}
 H\left(\pi_{L,N} \mid \bar{\nu}_{\varphi_L,L}^{\otimes L}\right)
&=
 \sum_{\eta\in \Omega_{L,N}} \pi_{L,N} [\eta] \log \left( \frac{\pi_{L,N} [\eta]}{\bar{\nu}_{\varphi_L,L}^{\otimes L}[\eta]} \right) 
=
\log \left( \frac{z_L(\varphi_L)^L }{Z_{L,N}\varphi_L^{N}} \right) \\
&=
-\log  \bar{\nu}_{\varphi_L,L}^{\otimes L}\left[ \sum_{x=1}^{L} \eta_x =N \right]\, .
\end{split}
\end{align*}
Then, by subadditivity of the relative entropy we have for marginals
\begin{equation}\label{eq_marg_rel_entropy}
H\left(\pi_{L,N}[\eta_x\in \cdot\,] \mid \bar{\nu}_{\varphi_L,L}\right)
\leq - \frac{1}{L} \log  \bar{\nu}_{\varphi_L,L}^{\otimes L}\left[ \sum_{x=1}^{L} \eta_x =N \right]\, .
\end{equation}
We estimate the right-hand side using the local limit theorem in Lemma \ref{lem_llt} with the specific choices of
\begin{align*}
    a_L:= L\, R_L(\varphi_L) \quad \text{and} \quad b_L:= \sqrt{L\, \sigma_L^2(\varphi_L)}\, .
\end{align*}
It follows from Lemma \ref{lem_check_lindeberg} that 
$$
\sum_{x=1}^L \zeta_{x,L} = \frac{\sum_{x=1}^L \eta_{x,L}-a_L}{b_L}\weakconv \mathcal{N}(0,1)\,.
$$
Moreover, $b_L$ diverges in the large $L$ limit because $\lim_{L\to\infty}\sigma_L^2(\varphi_L)=\sigma^2(\varphi)>0$.
Also, 
\begin{align*}
\limsup_{L\to\infty} \frac{b_L^2}{Q_L}
=
    \limsup_{L\to\infty} \frac{ \sigma_L^2(\varphi_L)}{q(\bar{\nu}_{\varphi_L,L})}\leq\frac{\sigma^2(\varphi)}{\sup_{n\in\mathbb{N}_0}[w(n)\wedge w(n+1)]}\,,
\end{align*}
where the final inequality follows by dominated convergence.
The right hand side is finite by assumption \ref{aspt_B1_uniform_conv}.
Therefore, we may apply the local limit theorem stated in Lemma \ref{lem_llt} which yields
$$
\bar{\nu}_{\varphi_L,L}^{\otimes L}\left[ \sum_{x=1}^{L} \eta_x =N \right]
= O\left(L^{-1/2}\right)\,,
$$
and coming back to \eqref{eq_marg_rel_entropy}
$$
H\left(\pi_{L,N}[\eta_x\in \cdot\,] \mid \bar{\nu}_{\varphi_L,L}\right)
\sim
\frac{\log L}{L}\to 0\quad\mbox{as }L\to\infty\, .
$$
With Pinsker's inequality (see e.g. \cite[Lemma 6.2]{Gray11})
this implies for the total variation distance 
\begin{align}
    \label{eq_TVconv}
    d_{TV} \left(\pi_{L,N}[\eta_x \in \cdot\,], \bar{\nu}_{\varphi_L,L}\right) \to 0\, .
\end{align}
Also, by \eqref{eq_moment_limits} and uniform convergence of the weights, $\bar{\nu}_{\varphi_L,L}$ converges weakly to $\bar{\nu}_\varphi$, which together with \eqref{eq_TVconv} implies
$
\pi_{L,N}[\eta_x \in \cdot\,]\weakconv \bar\nu_{\Phi (\rho )}
$ 
as $N/L\to\rho\leq \rho_c$.

Finally, we conclude the super-critical case $\rho >\rho_c$ using a large deviation estimate. Now let $\varphi_L$ be a sequence converging to $1$ and satisfying \eqref{eq_choice_phi} and \eqref{eq_moment_limits} so that $R_L(\varphi_L) \to \rho_c$.
Then
\begin{equation*}
\begin{split}
-\frac{1}{L} \log  \bar{\nu}_{\varphi_L,L}^{\otimes L}\left[ \sum_{x=1}^{L} \eta_x =N \right]
&\leq 
-\frac{1}{L} \log  \bar{\nu}_{\varphi_L,L}^{\otimes L\setminus \{1\}}\left[ \sum_{x=2}^{L} \eta_x =N- \lfloor(\rho-\rho_c)L\rfloor  \right] \\
&\qquad-\frac{1}{L} \log  \bar{\nu}_{\varphi_L,L}\left[  \eta_1 = \lfloor(\rho-\rho_c)L\rfloor \right]\, ,
\end{split} 
\end{equation*}
where the first term on the r.h.s.~converges to zero by the local central limit theorem, since $(N- \lfloor(\rho-\rho_c)L\rfloor )/L \to \rho_c$. 
For the second term,
\begin{align*}
&-\frac{1}{L} \log  \bar{\nu}_{\varphi_L,L}\left[  \eta_1 = \lfloor(\rho-\rho_c)L\rfloor \right] \\
&\qquad =  
-\frac{1}{L} \log  w_L(\lfloor(\rho-\rho_c)L\rfloor ) 
- (\rho-\rho_c) \log  \varphi_L 
+\frac{1}{L} \log  z_{L}(\varphi_L)\to 0 \, ,
\end{align*}
where convergence follows from the sub-exponential assumption \eqref{eq_sub_exp_tails}, and  since $\varphi_L \to 1$ and $z_{L}(\varphi_L)\to z(1)\in (0,\infty )$. 
It follows that $H\left(\pi_{L,N}[\eta_x\in \cdot\, ] \mid \bar{\nu}_{\varphi_L,L}\right)$ vanishes and, for the same reason as in the sub-critical case, $\pi_{L,N}[\eta_x\in \cdot\, ]$ converges weakly to $w=\bar\nu_1$.

Finally, to establish weak convergence of finite dimensional marginals; fix $n_1,\ldots,n_m\in\mathbb{N}_0$ and $x_1,\ldots,x_m$ distinct indices, then 
\begin{align*}
\pi_{L,N}&[\eta_{x_1}=n_1,\eta_{x_2}=n_2]\\
&=
\frac{w_L(n_1)Z_{L-1,N-n_1}}{Z_{L,N}}
\frac{1}{Z_{L-1,N-n_1}}
w_L(n_2)
\sum_{\xi \in \Omega_{L-2,N-n_1-n_2}} \left(\prod_{z} w_L(\xi_z)\right)\\
&=\pi_{L,N}[\eta_{x_1}=n_1]
\,\pi_{L-1,N-n_1}[\eta_{x_2}=n_2]\,.
\end{align*}
This identity immediately generalises to
\begin{align*}
\pi_{L,N}[\eta_{x_1}=n_1,\ldots,\eta_{x_m}=n_m]
= \prod_{j=1}^m 
\pi_{L-j+1,N-\sum_{k=1}^{j-1}n_j}[\eta_{x_j} = n_j]\,,
\end{align*}
which, by taking the thermodynamic limit on the right hand side, completes the proof.

It is possible to drop the assumption that the variance of the limiting weights is positive in the case $w(n)=\mathds{1}\{n=0\}$.
If we fix $\rho > 0$, we can use the same argument as in the super-critical case above and put all particles on a single site. 
In this case the right hand side of \eqref{eq_marg_rel_entropy} vanishes since
\begin{align*}
    - \frac{1}{L} \log  \bar{\nu}_{\varphi_L,L}^{\otimes L}\left[ \sum_{x=1}^{L} \eta_x =0 \right]
    = -\log  \nu_{\varphi_L,L} \left[ \eta_x =0 \right]\to \log w(0)=0\, ,
\end{align*}
as $L\to\infty$. In this case we do not use the local central limit theorem.
Otherwise, the proof remains unchanged.
\end{proof}



\medskip

\bibliography{cites}
\bibliographystyle{alpha}


\end{document}